\documentclass[journal,10pt]{IEEEtran}

\usepackage[utf8]{inputenc}

\usepackage{amsmath,amssymb,amsthm,mathtools,nccmath}
\usepackage{bm}
\usepackage{enumitem}
\usepackage{graphicx}

\usepackage{algorithmic}
\newcounter{algorithm}
\usepackage{xcolor}
\usepackage{microtype}
\usepackage{hyperref}

\hypersetup{hidelinks}
\usepackage{comment}

\theoremstyle{plain}
\newtheorem{theorem}{Theorem}
\newtheorem{lemma}{Lemma}

\newtheorem{corollary}{Corollary}
\theoremstyle{definition}
\newtheorem{assumption}{Assumption}
\newtheorem{definition}{Definition}
\theoremstyle{remark}
\newtheorem{remark}{Remark}

\newcommand{\ud}{\mathrm{d}}
\newcommand{\R}{\mathbb{R}}
\newcommand{\E}{\mathbb{E}}
\newcommand{\Prob}{\mathbb{P}}
\newcommand{\Ball}{\mathbb{B}}
\newcommand{\Sphere}{\mathbb{S}}
\newcommand{\Unif}{\mathrm{Unif}}
\newcommand{\vol}{\mathrm{vol}}
\newcommand{\inner}[2]{\langle #1,\, #2\rangle}
\newcommand{\norm}[1]{\left\lVert #1\right\rVert}
\newcommand{\xbar}{\bar{x}}
\newcommand{\gbar}{\bar{g}}
\newcommand{\bx}{\mathbf{x}}
\newcommand{\bg}{\mathbf{g}}
\newcommand{\bP}{\mathbf{P}}
\newcommand{\xperp}{\bx_{\perp}}
\newcommand{\one}{\mathbf{1}}
\newcommand{\Reg}{\mathcal{R}}
\newcommand{\Filt}{\mathcal{F}}
\newcommand{\grad}{\nabla}
\newcommand{\Dt}{\widetilde{D}}
\newcommand{\mut}{\widetilde{\mu}}
\newcommand{\gat}{\widetilde{\gamma}}
\DeclareMathOperator*{\argmin}{arg\,min}

\DeclarePairedDelimiterX\Econd[2]{[}{]}{#1 \mkern2mu\delimsize\vert\mkern2mu #2}
\DeclarePairedDelimiterX\Pcond[2]{(}{)}{#1 \mkern2mu\delimsize\vert\mkern2mu #2}

\newcommand{\Lcal}{\mathcal{L}}
\definecolor{revisionblue}{RGB}{0,76,153}
\colorlet{revisionblue}{black}

\newenvironment{revision}{\color{revisionblue}}{}

\title{Differentially Private Distributed Optimization via the Intrinsic Randomness of Zeroth-Order Estimators: Sampling Geometry, Convergence, and Privacy Guarantees}
\author{Jiajun Wang and Yujie Tang%
\thanks{This work was supported by the National Natural Science Foundation of China through grants 72431001 and 72301008.
The authors are with the Department of Control Science and Systems Engineering, Peking University. Email: {\tt\small wangjiajun@stu.pku.edu.cn}, {\tt\small yujietang@pku.edu.cn}}}
\begin{document}
\maketitle

\begin{abstract}
Distributed zeroth-order optimization has gained increasing attention in applications where exact gradients are computationally expensive or unavailable. However, the iterative exchange of local decision variables makes distributed networks vulnerable to inference attacks. While differential privacy provides a rigorous framework for privacy protection, conventional approaches generally require additional random noise, which inevitably affects optimization accuracy. In this paper, we investigate whether the randomness inherent in zeroth-order estimation can itself provide privacy protection. We develop a distributed zeroth-order gradient descent method that integrates consensus updates with a ball-uniform two-point gradient estimator. Our theoretical analysis shows that our algorithm attains standard convergence guarantees in both online and offline optimization settings. We then characterize the relevant probability distributions of the associated privacy mechanism, and establish differential privacy guarantees for both the online and offline optimization settings. These results reveal the potential of intrinsic zeroth-order randomness for privacy protection.
\end{abstract}

\section{Introduction}
\label{sec:intro}

Distributed optimization considers a network of agents that cooperatively
minimizes an aggregate objective while each agent has access only to its own
local objective and communicates only with neighboring agents. Early
consensus and subgradient methods established the basic paradigm of combining
local descent with weighted averaging
\cite{nedic2009distributed,nedic2010constrained}. Subsequent studies
characterized the behavior of decentralized gradient descent
\cite{yuan2016convergence} and developed fixed-step methods, including
EXTRA and gradient-tracking algorithms
\cite{shi2015extra,nedic2017diging,qu2018harnessing}. These developments have
made distributed optimization an important computational framework for
networked control, estimation, and learning. When the local objectives change
over time and arrive sequentially, the agents must make decisions before future objectives are known,
leading to distributed online optimization. Its performance is commonly
measured by regret relative to the best fixed decision
\cite{hosseini2013online}. In both the offline and online settings, however,
most classical distributed algorithms presume that every agent can obtain
a gradient or subgradient of its local objective.

This first-order requirement can be restrictive when local objectives are
available only through function evaluations. Zeroth-order optimization
addresses this limitation by querying an objective at randomly perturbed
points and using the resulting function differences to construct a
gradient estimator. Random smoothing and gradient estimators
form the basis of zeroth-order and bandit optimization
\cite{flaxman2005online,duchi2015optimal,shamir2017optimal,
nesterov2017random}. Combining zeroth-order estimation with distributed optimization introduces
additional difficulties. Each
agent constructs an estimator from its own random query direction, while the
resulting estimation error propagates through the communication network. To address these difficulties, distributed
zeroth-order methods have been studied for constrained convex optimization
\cite{yuan2016zeroth}, random networks
\cite{sahu2018distributed}, nonconvex multi-agent problems
\cite{hajinezhad2019zone}, consensus and gradient-tracking methods
\cite{tang2021distributed}, and stochastic distributed optimization
\cite{yi2022zeroth,yuan2024distributed}.

A separate challenge arises from the communication process of distributed optimization.
Although the objectives and raw data remain local, the sequence of
transmitted states depends on the agents' private objectives and
may reveal information about them. Differential privacy (DP) offers a
mathematically rigorous framework that quantifies and limits
the privacy risk associated with a prescribed change of the problem instance
\cite{dwork2006calibrating,dwork2014algorithmic}. In distributed
optimization, privacy is commonly enforced by adding calibrated randomness
to communicated states, local objectives, gradients, or update directions
\cite{huang2015differentially,nozari2016differentially,
ding2021differentially,wang2024tailoring}. Differential privacy and
zeroth-order feedback have also been combined in decentralized optimization
by perturbing both the agents' states and their gradient-tracking variables
\cite{zhao2024privateonline}. While such mechanisms provide rigorous
protection, the additional noise enters the optimization recursions, and
repeated releases lead to an accumulated privacy budget. This creates a
trade-off between privacy and optimization accuracy.

The randomness of a zeroth-order update offers a possible
way to reduce the cost of explicit privatization. Recent zeroth-order methods for private
large-model fine-tuning replace the release of a high-dimensional private
gradient by a clipped and noised scalar directional statistic, thereby
avoiding backpropagation and per-example-gradient storage
\cite{zhang2024dpzero,tang2024private}. Privacy accounting can be further
improved when intermediate optimization states remain hidden
\cite{chien2025privacy}. These developments raise a more fundamental
question: Can the randomness in a zeroth-order estimator
itself protect the released updates, without an independent privacy
perturbation? Gratton et al.~\cite{gratton2022privacy} provided an affirmative indication
in a distributed zeroth-order ADMM framework. Their privacy argument
approximates a batch-averaged zeroth-order perturbation by a Gaussian
distribution through the central limit theorem. This requirement is not
innocuous: Obtaining an accurate Gaussian approximation requires a
sufficiently large number $K$ of random-direction samples, which increases
the function-query cost while reducing the covariance of the intrinsic
perturbation by a factor $1/K$. Thus, the condition used to justify the
Gaussian approximation simultaneously weakens the randomness available for
privacy, and the finite-$K$ approximation error is not included in the
privacy guarantee. On the other hand, Gupta et
al.~\cite{gupta2025inherent} showed that zeroth-order randomness does not
automatically provide a universal privacy guarantee. Their lower bounds use
worst-case adjacent objectives for which a zero-preserving oracle produces a
deterministic zero update under one instance and an atomless update under the
other. These results rule out unconditional inherent privacy over an
unrestricted neighboring class, but do not exclude conditional privacy when
the gradients generated by admissible adjacent objectives are nondegenerate
and relatively close.

Motivated by this gap, this paper investigates differential privacy guarantees that stem from the inherent randomness of zeroth-order gradient estimators. Our study is rooted in the support geometry of a
ball-uniform two-point estimator. We develop a consensus-based distributed
zeroth-order method and study its online regret and offline convergence. For
privacy, we analyze the associated limit-form algorithm using the directional-derivative estimator
under online and offline adjacency. The main
contributions are summarized as follows.

\begin{itemize}
    \item
    We identify support geometry of the random direction as a critical feature that characterizes whether the randomness of a zeroth-order estimator can intrinsically give rise to differential privacy. We then show that sampling uniformly from the unit ball, rather than from the unit sphere, yields a limit-form distribution suitable
    for nontrivial privacy. For the ball-uniform
    directional-derivative estimator, we derive its density and support, reduce the likelihood-ratio calculation for adjacent gradients to a
    one-dimensional Cauchy variable, and establish an explicit single-round \((\epsilon,\delta)\)-differential privacy guarantee under a gradient-discrepancy condition.

    \item
    We establish optimization guarantees for our proposed algorithm using the ball-uniform estimator. We first establish fundamental statistical properties of the ball-uniform estimator and analytical properties of the associated smoothed function. Then, for online convex and Lipschitz losses, the
    proposed method achieves $O(\sqrt{T})$ expected static network regret. For
    offline optimization with strongly convex and smooth local objectives, we prove an
    $O(1/t)$ expected last-iterate bound on the squared distance to the minimizer, together with an $O(1/t^2)$ mean-square consensus error bound.

    \item
    We extend the single-round privacy result to the complete broadcast trajectory where the algorithm employs the directional-derivative estimator. Under online adjacency, only one conditional round differs, so the
    resulting privacy bound is independent of the horizon $T$ and the
    number of agents $N$, whenever the single-round value satisfies
    $\delta_{\rm step}<1$. Under offline adjacency, the affected
    local objective differs at every round, and we apply adaptive advanced
    composition, which can improve
    the dependence on $T$ for suitable privacy parameters, at the cost of an
    additional slack $\widehat\delta$.
\end{itemize}

\paragraph*{Notation} For each positive integer $N$, we let $[N]$ denote the set $\{1,\ldots,N\}$. We use $\|x\|$ to denote the standard $\ell_2$ norm of a real vector $x$, use $\|x\|_1$ to denote its $\ell_1$ norm, and use $\|A\|_2$ to denote the spectral norm of a matrix $A$. Given a list of vectors $x_1,\ldots,x_N$, we use $\mathrm{col}(x_i)_{i\in[N]}$ to denote the vector formed by concatenating $x_1,\ldots,x_N$ together. The $j$th standard basis vector is denoted by $e_j$. We let $\Sphere^{d-1} \coloneq \{w \in \R^d : \norm{w} = 1\}$ denote the unit sphere
in $\R^d$, and let $\Ball^d \coloneq \{w \in \R^d : \norm{w} \le 1\}$ denote the closed unit
ball.

\section{Problem Formulation and Preliminaries}
\label{sec:preliminaries}

\subsection{Problem Formulation}
\label{sec:prob}

Let the set of agents be $[N] = \{1, \dots, N\}$. The agents
communicate over a network whose topology is specified by an undirected graph
$\mathcal{G} = ([N], \mathcal{E})$, with each edge $(i,j) \in \mathcal{E}$
representing a bidirectional communication link. Agent~$i$ exchanges messages only
with its neighbors in the set $\mathcal{N}_i \coloneq \{j : (i,j) \in \mathcal{E}\}$. We adopt a weight matrix $P\in\R^{N\times N}$ that will be used in our algorithm design.

\begin{assumption}\label{as:network}
The graph $\mathcal{G}$ is connected. The matrix
$P\in\R^{N\times N}$ is symmetric and doubly
stochastic, and satisfies $P_{ij}=0$ whenever
$i\ne j$ and $(i,j)\notin\mathcal E$. Moreover,
\[
  \sigma\coloneq \left\lVert P-\frac1N\one\one^{\!\top}\right\rVert_2<1.
\]
\end{assumption}

In this article, we investigate two problem settings that share the network under
Assumption~\ref{as:network}: the online setting and the offline setting.
In both settings, each agent
accesses its private local objective only through a zeroth-order oracle, i.e., it
can query the values of the objective at finitely many points but does not
observe its gradient.

\subsubsection{Online Setting}
In the online setting, there is a discrete time horizon $\mathcal T\coloneq \{0,1,\ldots,T-1\}$ where $T$ is a given positive integer; we also denote $\mathcal T_+ \coloneq \{0,1,\ldots,T\}$.
At the start of round~$t\in\mathcal T$,
each agent~$i$ holds a state $x_i(t)\in\R^d$. Then, agent~$i$ is presented with a private local loss
function $f_{i,t} : \R^d \to \R$, drawn from an admissible class $\Lcal$ of
convex Lipschitz losses (see Assumption~\ref{as:reg} below). We refer to the loss
collection
\begin{equation}\label{eq:problem}
  F \;:=\; \{f_{i,t}\}_{i \in [N],\, t \in \mathcal T} \;\in\; \Lcal^{N\times T}
\end{equation}
as a \emph{problem instance}. The instance is determined before round~$0$
and does not depend on the algorithm's execution.

\begin{assumption}\label{as:reg}
The online loss class $\Lcal$ consists of convex, continuously differentiable
functions $f:\R^d\to\R$ that are globally $L_0$-Lipschitz:
$|f(x)-f(y)|\le L_0\norm{x-y}$ for all $x,y\in\R^d$, where $L_0>0$.  Every online loss
$f_{i,t}$ belongs to $\Lcal$.
\end{assumption}

Assuming that the cumulative loss admits a minimizer, we fix a comparator
\begin{equation}\label{eq:online-comparator}
x^\star \in \argmin_{x \in \R^d}
\sum_{t=0}^{T-1}\sum_{i=1}^{N} f_{i,t}(x),
\end{equation}
and measure performance by the static network regret for
distributed online learning, which evaluates the
decisions of a single agent on the losses of the whole network \cite{yan2013distributed}.

\begin{definition}\label{def:regret}
For any agent $j \in [N]$, the regret incurred by the iterate sequence
$\{x_j(t)\}_t$ is
\begin{equation}\label{eq:regret-def}
  \Reg_T(j) \;\coloneq\;
  \sum_{t=0}^{T-1} \sum_{i=1}^{N}
    \bigl[\,f_{i,t}(x_j(t)) - f_{i,t}(x^\star)\,\bigr].
\end{equation}
\end{definition}

\subsubsection{Offline Setting}
In the offline setting, each agent~$i$ holds a fixed private local objective $f_i:\R^d\to\R$ that
is used at every round of our iterative algorithm, and cooperatively solves
\begin{equation}\label{eq:offline-model}
  \min_{x\in\R^d}\ f(x)\coloneq\frac1N\sum_{i=1}^N f_i(x).
\end{equation}
The problem instance is $F^{\mathrm{off}}\coloneq(f_1,\ldots,f_N)$. Assuming that $f$ has a minimizer $x^\star$, we measure performance at each round
$t$ by the expected optimality gap of the network average
$\xbar(t)\coloneq\frac1N\sum_{i=1}^N x_i(t)$ compared to $x^\star$ and by the consensus
error $\frac1N\sum_{i=1}^N\E\norm{x_i(t)-\xbar(t)}^2$.

In the following, we shall use a unified notation $f_{i,t}$ to denote the local objective of agent $i$ for both the online and the offline settings, with $f_{i,t}=f_i,\,\forall t$ in the offline setting.

\subsubsection{Privacy Model}
We consider privacy threats caused by an external curious adversary
who observes the entire sequence of broadcast states. Specifically, the
adversary observes and records the trajectory
\begin{equation}
\label{eq:obs}
  \Pi_T \;\coloneq\; \bigl(x_i(t)\bigr)_{i \in [N],\, t \in \mathcal T_+}. 
\end{equation}
The deterministic initialization, the weight matrix $P$, and the schedules $\{\eta_t\}$ and
$\{\beta_t\}$ are public, fixed, and independent of the objectives. The
queried function
values are internal and are not revealed; the adversary observes only the
broadcast trajectory $\Pi_T$. This type of privacy model has been prevalent in existing literature on differentially private distributed optimization~\cite{huang2015differentially,nozari2016differentially,ding2021differentially,wang2024tailoring}.

The goals are twofold: (i) to develop distributed zeroth-order algorithms for the online and offline settings with corresponding regret convergence rate guarantees; and (ii) to establish
$(\epsilon,\delta)$-differential privacy in the settings above.

\subsection{Preliminaries on Zeroth-Order Optimization}
\label{sec:prelim_zo}

In zeroth-order optimization, the decision-maker can evaluate a loss function but
cannot access its derivatives. The gradient used by first-order methods is
therefore replaced by a randomized estimate: The query
point is perturbed along a random direction, the resulting pair of function
values is compared, and the direction is rescaled by the observed
difference.

Given a function $h : \R^d \to \R$, a smoothing radius $\beta > 0$,
and a random direction $u$ sampled from a distribution $\mathcal{D}$
on $\R^d$, we consider the following zeroth-order gradient estimator:
\begin{equation}\label{eq:est-prelim}
  G_h(x;\beta,u)
  \;\coloneq \;
  c_{\mathcal{D}}\cdot \frac{h(x+\beta u) - h(x-\beta u)}{2\beta}\, u ,
\end{equation}
where $c_{\mathcal{D}}>0$ is a normalization constant matched to the distribution $\mathcal{D}$. Two classical choices of $\mathcal{D}$ are the
sphere-uniform direction $u\sim\Unif(\Sphere^{d-1})$ with
$c_{\mathcal{D}}=d$, and the Gaussian direction $u\sim\mathcal N(0,I_d)$ with $c_{\mathcal{D}}=1$. In each case, the constant $c_{\mathcal{D}}$ is calibrated so
that $\E\bigl[G_h(x;\beta,u)\bigr]$ equals the gradient of a smoothed
version of $h$ \cite{shamir2017optimal,nesterov2017random}.

To expedite our theoretical analysis of the intrinsic privacy guarantees of zeroth-order gradient estimators, we shall introduce a limit version of~\eqref{eq:est-prelim}. For a differentiable $h$, as $\beta\downarrow 0$, we have
$\tfrac{1}{2\beta}\bigl(h(x+\beta u)-h(x-\beta u)\bigr)\to
u^{\!\top}\grad h(x)$, so that the estimator reduces to
\[
  G_h(x;\beta,u)
  \;\longrightarrow\;
  c_{\mathcal{D}}\,\bigl(u^{\!\top}\grad h(x)\bigr)\,u .
\]
We call $c_{\mathcal{D}}\,\bigl(u^{\!\top}\grad h(x)\bigr)\,u$ the \emph{directional-derivative estimator}, which has also been studied in \cite{nesterov2017random}. This limit version can indeed facilitate our differential privacy analysis, as the densities of the estimator and the associated likelihood ratio will then have closed-form expressions (Lemmas~\ref{lem:density} and~\ref{lem:cauchy}), serving as foundations for establishing differential privacy guarantees.
Strictly speaking, however, it is a different estimator from the finite-difference
version~\eqref{eq:est-prelim} with $\beta>0$.  Accordingly, the privacy theorems in
Section~\ref{sec:privacy} are stated only for this directional-derivative estimator. Nevertheless, we expect that the results in Section~\ref{sec:privacy} will provide novel and nontrivial insights on how the intrinsic randomness of zeroth-order gradient estimators can lead to differential privacy guarantees.
The convergence theorems, on the other hand, concern the implementable finite-difference estimator~\eqref{eq:est-prelim}.

For optimization purposes, the choice of the distribution $\mathcal{D}$ is largely
immaterial: Each of the common choices yields the standard convergence
guarantees of two-point zeroth-order methods~\cite{nesterov2017random,shamir2017optimal}. For differential privacy, however, it is decisive: The choice of
$\mathcal{D}$ determines the geometry of the estimator's distribution, and with it whether differential privacy is achievable at all.
This choice is the central design decision of our algorithm, and we devote
Section~\ref{sec:design} to it.

\subsection{Preliminaries on Differential Privacy}
\label{sec:prelim_privacy}

We consider the case where the problem instance itself is the sensitive object: The online losses
$f_{i,t}$, or the offline objectives $f_i$, encode the agents' private information. Accordingly, we define adjacency at the
level of problem instances.

\begin{definition}[Online adjacency]\label{def:adj}
Two problem instances $F = \{f_{i,t}\}$ and $F' = \{f'_{i,t}\}$ in
$\Lcal^{N\times T}$ are \emph{adjacent}, written $F \sim F'$, if there
exists a unique index pair $(i_0, t_0) \in [N] \times \mathcal T$ such that
$f_{i_0, t_0} \neq f'_{i_0, t_0}$, while $f_{i,t} = f'_{i,t}$ for every
other $(i, t) \in [N] \times \mathcal T$.
\end{definition}

\begin{definition}[Offline adjacency]
\label{def:offline-adj}
Two problem instances
\[
  F^{\mathrm{off}}=(f_1,\ldots,f_N),
  \qquad
  (F')^{\mathrm{off}}=(f'_1,\ldots,f'_N)
\]
are \emph{adjacent} if there is a unique $i_0\in[N]$ such that
$f_{i_0}\ne f'_{i_0}$, while $f_i=f'_i$ for every $i\ne i_0$. Thus the
same affected local function is used at every round.
\end{definition}

For a problem instance $F$, a specified randomized update rule and
release protocol induce an output sequence $\mathcal M(F)$. Its privacy
is measured by the following standard definition.

\begin{definition}[$(\epsilon,\delta)$-differential privacy]\label{def:dp}
Let $\sim$ be a specified adjacency relation on an input space. For
$\epsilon\ge0$ and $\delta\in[0,1]$, a randomized mechanism $\mathcal M$
with values in a measurable output space is
\emph{$(\epsilon,\delta)$-differentially private} if, for every adjacent pair
$F\sim F'$ and every measurable output event $\mathcal S$,
\begin{equation}\label{eq:dp}
  \Prob\bigl[\mathcal{M}(F) \in \mathcal{S}\bigr]
  \;\le\;
  e^{\epsilon}\,
  \Prob\bigl[\mathcal{M}(F') \in \mathcal{S}\bigr]
  \;+\;\delta .
\end{equation}
\end{definition}

\section{Ball-Uniform Zeroth-Order Distributed Algorithm}
\label{sec:algorithm}

\subsection{Design of the Random Directions}
\label{sec:design}

Before stating the algorithm, we explain how its privacy mechanism is
designed. The guiding question is the following: \emph{Can the randomness 
from a zeroth-order estimator be shaped into a differential privacy guarantee,
without injecting any additional noise?} By answering it, one can fix, step by step,
each critical design choice of our algorithm.
 
Consider a classical distributed gradient descent iteration in which agent $i$
updates
$x_i(t+1)=\sum_{j} P_{ij}\,x_j(t)-\eta_t\,g_i(t)$
and broadcasts the result to its neighbors, where $g_i(t)$ is some randomized gradient
estimator. Conditioned on the history up to the start of round~$t$, the consensus term $\sum_j P_{ij}x_j(t)$ becomes known, and the broadcast $x_i(t+1)$ is the image of $g_i(t)$ under
a fixed invertible map
$y\mapsto\sum_j P_{ij}x_j(t)-\eta_t\,y$ whenever $\eta_t>0$, as assumed
throughout.
Invertible affine maps leave likelihood ratios unchanged, so the
likelihood ratio between the distributions of $x_i(t+1)$ under two adjacent
problems, which is a fundamental quantity in differential privacy analysis, agrees with that between the distributions of the corresponding gradient estimator $g_i(t)$. The design question therefore becomes: Which
zeroth-order estimator has an output distribution suitable for differential
privacy?
 
\begin{figure*}[t]
  \centering
  \includegraphics[width=\linewidth]{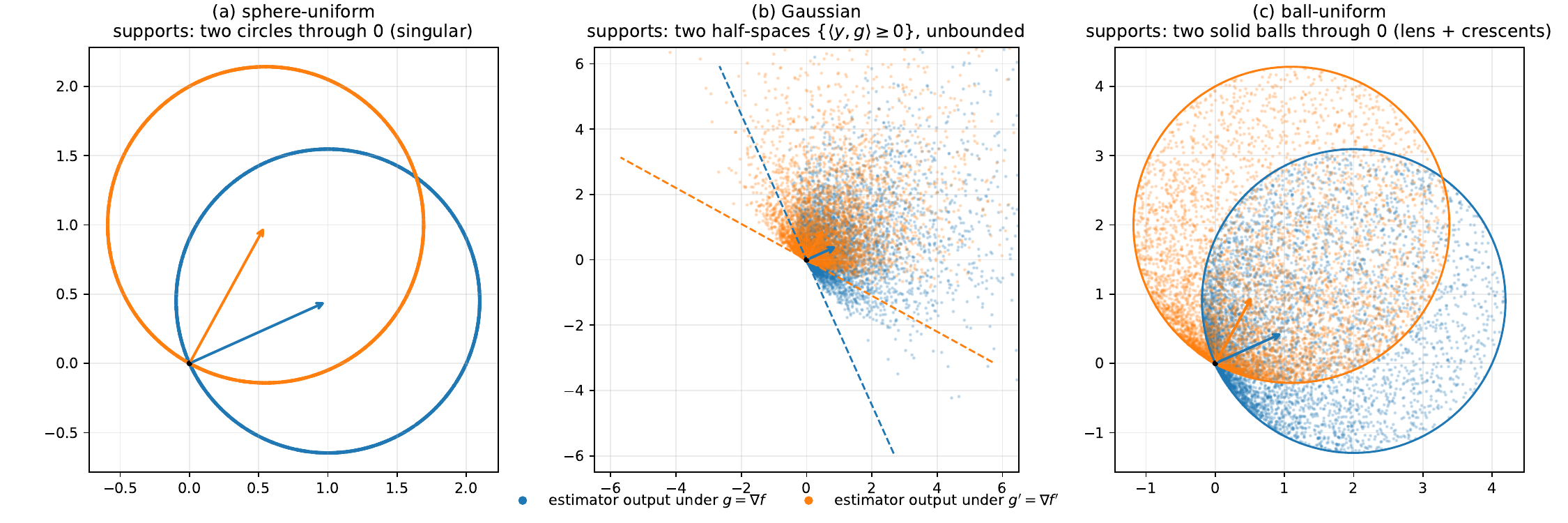}
  \caption{Probability distributions of the (limit-form) two-point estimator in
  $d=2$, under two adjacent problems with gradients $g$ (blue) and $g'$
  (orange); $6{,}000$ samples each.
  (a)~Sphere-uniform sampling: Each distribution concentrates on a circle through
  the origin; the two circles
  share only a null set, the distributions are mutually singular, and differential
  privacy requires $\delta\ge1$ for this distinct-gradient pair.
  (b)~Gaussian sampling: The distributions have densities, with the
  supports being the half-spaces $\{\inner{y}{g}\ge0\}$ (dashed boundaries).
  (c)~Ball-uniform sampling: The supports are the solid balls
  $\overline{B}\bigl(\tfrac{c_{\mathcal{D}}}{2}g,\tfrac{c_{\mathcal{D}}}{2}\!\norm{g}\bigr)$ through the
  origin. In the displayed non-nested case, the lens-shaped overlap
  carries both densities, while the two crescents are each covered by only
  one distribution; this nonoverlap motivates the high-risk region
  $D_\epsilon(g,g')$.}
  \label{fig:supports}
\end{figure*}

The classical two-point estimator draws its random direction $u$ uniformly from the
unit sphere $\Sphere^{d-1}$. For a differentiable function $h$, when the smoothing radius $\beta$ is sufficiently small, the gradient estimator~\eqref{eq:est-prelim} can be approximated by the directional-derivative estimator
\[
  y \;=\; d\bigl(u^{\!\top}\grad h(x)\bigr)\,u ,
\]
as discussed in Section~\ref{sec:prelim_zo}. We then ask how the above $y$ responds when the function $h$ is replaced by
$h'$, and with it the gradient $g=\grad h(x)$ by $g'=\grad h'(x)$. The geometry of the support of the gradient estimator
already settles the question: 
Straightforward calculation shows that $\|y-\frac{d}{2}g\|^2=(\frac{d}{2}\|g\|)^2$, implying that $y$ lies on the $(d-1)$-dimensional sphere passing through the origin with center $\tfrac d2 g$ and diameter $d\norm{g}$. Thus, the gradient estimator $y$ is concentrated on a set of zero Lebesgue measure
(see Figure~\ref{fig:supports}(a)). When $h$ is replaced by $h'$ with a different gradient, the estimator $y'=d\bigl(u^{\!\top}\grad h'(x)\bigr)\,u$ lives
on a \emph{different} such sphere, and two distinct spheres through the
origin intersect in a set that carries no probability mass under either distribution.
The distributions of $y$ and $y'$ are thus mutually singular, i.e., there is an event that
has probability~1 under one problem and probability~0 under the other. Consequently, the differential privacy inequality \eqref{eq:dp} can
only hold with $\delta\ge1$ for any $\epsilon>0$. We see that the
gradient estimator with random direction $u\sim\Unif(\Sphere^{d-1})$ cannot protect two distinct gradients at a
nontrivial privacy level unless a further randomization mechanism (such as adding independent Laplace-distributed noise) is incorporated.
 
Fortunately, the above impossibility argument also suggests a solution: 
Unlike the uniform distribution on the unit sphere $\Sphere^{d-1}$, if the random direction $u$ instead follows a distribution with a density on $\R^d$, then one may expect that the resulting distribution of $y$ will admit finite likelihood ratio under adjacent problems. Two natural candidates are
the Gaussian direction $u\sim\mathcal N(0,I_d)$ and the ball-uniform
direction $u\sim\Unif(\Ball^d)$, both yielding estimators of the form
$y\propto(u^{\!\top}g)\,u$.

Both candidates are viable in principle. For Gaussian directions, however, 
the associated likelihood ratios are unwieldy. For ball-uniform directions, the density of $(u^{\!\top}g)\,u$ admits a simple, elementary closed-form expression (derived
in Section~\ref{sec:privacy}), and the resulting density ratios reduce to
tractable one-dimensional quantities. We therefore build our algorithm on
the ball-uniform direction, i.e., we let the random direction $u$ follow $\Unif(\Ball^d)$ in the gradient estimator.

For $c\in\R^d$ and $r\ge0$, write
$\overline{B}(c,r):=\{y\in\R^d:\norm{y-c}\le r\}$.
The support of the ball-uniform estimator
$y=c_{\mathcal{D}}(u^{\!\top}g)\,u$ is the \emph{solid} ball through the origin
\[
  \mathrm{supp}(y)
  \;=\;
  \overline{B}\Bigl(\frac{c_{\mathcal{D}}}{2}\,g,\;\frac{c_{\mathcal{D}}}{2}\norm{g}\Bigr),
\]
as shown in Figure~\ref{fig:supports}(c). When the function $h$ is changed to a different $h'$, the ball remains
anchored at the origin but tilts and rescales with the gradient.  The supports may be non-nested, as in
Figure~\ref{fig:supports}(c), or nested when the gradients are positively
collinear.  In either case, at least one distribution assigns positive mass to a
region on which the other density vanishes. On such a region the likelihood
ratio is infinite, and so a pure $\epsilon$-guarantee is still impossible, indicating that some probability mass must be given
up and attributed to the additive slack $\delta$ in $(\epsilon,\delta)$-differential privacy. This dictates both the target and the proof strategy: We aim for
$(\epsilon,\delta)$-differential privacy, and we quantify the give-up mass
through the \emph{high-risk region}
\[
D_\epsilon(g,g')
  \coloneq \bigl\{y\in\R^d:\phi_g(y)>e^\epsilon\phi_{g'}(y)\bigr\},
\]
where $\phi_g$ and $\phi_{g'}$ denote the Lebesgue densities of the
(unnormalized) estimators $(u^{\!\top}g)u$ and $(u^{\!\top}g')u$, given in closed
form in Lemma~\ref{lem:density}. Since $y\mapsto c_{\mathcal{D}}y$ is an invertible
linear map, the likelihood ratio, and hence every quantity defined through
it below, is the same for the normalized estimators $c_{\mathcal{D}}(u^{\!\top}g)u$ and $c_{\mathcal{D}}(u^{\!\top}g')u$.
Section~\ref{sec:privacy}
computes the density in closed form, bounds the high-risk mass explicitly,
and shows that the resulting bound on $\delta$ tends to zero as the two gradients
approach one another.
 
We therefore adopt the ball-uniform direction $u\sim\Unif(\Ball^d)$ for the gradient estimator. For each $t$, agent $i$ draws an
independent direction $u_i(t)\sim\Unif(\Ball^d)$ and forms the ball-uniform
version of the gradient estimator \eqref{eq:est-prelim}:
\begin{equation}\label{eq:est-agent}
\begin{aligned}
 g_i(t):={}&\frac{d+2}{2\beta_t}
 \bigl[f_{i,t}\bigl(x_i(t)+\beta_t u_i(t)\bigr)\\
 &\qquad-f_{i,t}\bigl(x_i(t)-\beta_t u_i(t)\bigr)\bigr]u_i(t).
\end{aligned}
\end{equation}
Here we have plugged in the normalization constant $c_{\mathcal{D}} = d+2$ for the ball-uniform version; it offsets the second moment $\E\norm{u}^2 = d/(d+2)$ of the ball-uniform direction, and Lemma~\ref{lem:ball}(i) below will soon confirm that with this choice the estimator
is unbiased for a smoothed version of the loss.


The following lemma records the fundamental properties of the ball-uniform estimator and its associated smoothed function;
its proof is given in Appendix~\ref{app:estimator}.

\begin{lemma}\label{lem:ball}
Let $h:\R^d\to \R$ be continuous. For each $x\in\R^d$ and $\beta>0$, consider the gradient estimator
\[
G_h(x;\beta,u) = (d+2)\cdot\frac{h(x+\beta u)-h(x-\beta u)}{2\beta}\,u,
\]
where $u\sim\Unif(\Ball^d)$. Furthermore, define
\begin{equation}\label{eq:Phi}
  h^\beta(x)\coloneq (d+2)\int_0^1 r^{d+1}\,\hat h^{\beta r}(x)\,\ud r
\end{equation}
for all $x\in\R^d$, where for each $s>0$ we denote
\begin{equation}\label{eq:smooth}
  \hat h^{s}(x)\coloneq \E_{w\sim\Unif(\Ball^d)}\bigl[h(x+s w)\bigr].
\end{equation}
Let $u\sim\Unif(\Ball^d)$. Then the following results hold:
\begin{enumerate}[label=\textup{(\roman*)},leftmargin=2.4em,itemsep=1pt]
  \item $\E_{u}[G_h(x;\beta,u)]=\grad h^\beta(x)$.
  \item If $h$ is convex, then $h^\beta$ is also convex.
  \item If $h$ is $L_0$-Lipschitz, then $h^\beta$ is also $L_0$-Lipschitz, and
\begin{align*}
\bigl|h^\beta(x)-h(x)\bigr| & \le\beta L_0, \\
\E_{u}\!\left[\|G_h(x;\beta,u)\|^2\right]
& \le L_0^2(d+4).
\end{align*}
\item If $h$ is $L_1$-smooth, i.e., $\|\nabla h(x)-\nabla h(y)\|\leq L_1\|x-y\|,\,\forall x,y$, then $h^\beta$ is also $L_1$-smooth. Moreover, for all $x$, we have $\norm{\nabla h^\beta(x)-\nabla h(x)}\leq\beta L_1$ and 
\begin{align*}
& \E_u\!\left[
\norm{G_h(x;\beta,\! u)}^2
\right] \\
\le{} & 2(d\!+\!1)\|\nabla h(x)\|^2 + 2(d\!+\!2)\beta^2 L_1^2.
\end{align*}
\end{enumerate}
\end{lemma}

\subsection{The Algorithm}

Substituting the estimator~\eqref{eq:est-agent} into the classical decentralized
gradient descent framework yields our algorithm, \textsf{ZO-DGD},
stated in detail in Algorithm~\ref{alg:ozodgd}. It injects no additional noise for the purpose of privacy preservation.
While the design of Algorithm~\ref{alg:ozodgd}, after specifying the random direction $u$ in the gradient estimator, is straightforward (see, e.g., \cite{tang2021distributed,yi2022zeroth}), we still need to answer the following questions:
\begin{enumerate}
\item Can we still establish sublinear regret bounds or convergence rate guarantees for Algorithm~\ref{alg:ozodgd}, when we adopt the ball-uniform estimator?
\item Does Algorithm~\ref{alg:ozodgd} satisfy $(\epsilon,\delta)$-differential privacy? If so, how do $\epsilon$ and $\delta$ depend on the problem instances?
\end{enumerate}
In Section~\ref{sec:theory}, we provide theoretical convergence analysis of Algorithm~\ref{alg:ozodgd} with the finite-difference estimator~\eqref{eq:est-agent}, while
Section~\ref{sec:privacy} studies the variant where each $g_i(t)$ is replaced by the directional-derivative estimator and
obtains privacy guarantees from solely the intrinsic randomness.
 
\begin{figure}[t]
\refstepcounter{algorithm}
\label{alg:ozodgd}
\textbf{Algorithm \thealgorithm.} Zeroth-Order Decentralized Gradient
Descent (\textsf{ZO-DGD})
\par\smallskip\hrule\smallskip
\begin{algorithmic}[1]
\REQUIRE
  Initial states $\{x_i(0)\}_{i\in[N]}$, positive step-size schedule
  $\{\eta_t\}_{t\ge0}$, positive smoothing radius schedule
  $\{\beta_t\}_{t\ge0}$, weight matrix $P$.
\STATE \textbf{Initialize communication:} Each agent $i$ broadcasts
  $x_i(0)$ to its neighbors in $\mathcal{N}_i$ and receives the neighbors' initial states.
\FOR{$t = 0, 1, \dots, T-1$}
  \FOR{each agent $i \in [N]$ in parallel}
    \STATE Sample a private direction
      $u_i(t) \sim \Unif(\Ball^d)$.
    \STATE \textbf{Form gradient estimate:} query $f_{i,t}(\cdot)$ at
      $x_i(t)\pm\beta_t u_i(t)$ and form $g_i(t)$ according to
      \eqref{eq:est-agent}.
    \STATE \textbf{Update:}
      \[
        x_i(t+1) \;=\;
          \sum_{j=1}^{N} P_{ij}x_j(t)-\eta_t\,g_i(t) .
      \]
    \STATE \textbf{Communicate:} Broadcast $x_i(t+1)$ to neighbors in $\mathcal{N}_i$ and receive
      $\{x_j(t+1)\}_{j\in\mathcal N_i}$ for the next round.
  \ENDFOR
\ENDFOR
\end{algorithmic}
\smallskip\hrule
\end{figure}

\subsection{Matrix and Probability Notation}
 
For convenience of analysis, we collect the agents' iterates and gradient estimates into
stacked vectors
\[
\begin{aligned}
 \bx(t) \coloneq \mathrm{col}\bigl(x_i(t)\bigr)_{i\in[N]},\qquad
 \bg(t) \coloneq \mathrm{col}\bigl(g_i(t)\bigr)_{i\in[N]},
\end{aligned}
\]
and write $\bP := P\otimes I_d \in \R^{Nd\times Nd}$ for the Kronecker product of
the weight matrix $P$ and the identity matrix $I_d$. The averaged iterate and averaged gradient estimate are
\[
  \xbar(t) \;\coloneq\; \frac{1}{N}\sum_{i=1}^{N}x_i(t),
  \qquad
  \gbar(t) \;\coloneq\; \frac{1}{N}\sum_{i=1}^{N}g_i(t),
\]
and the consensus-error vector is
$\xperp(t) \coloneq \bx(t)-\one\otimes\xbar(t)$.
We also denote
\[
  J\coloneq \frac1N\one\one^{\!\top},\qquad
  M\coloneq I_N-J.
\]
Then $\norm{M}_2\le1$, $\norm{P-J}_2=\sigma$, and $MP=PM=P-J$.
All directions $\{u_i(t)\}_{i,t}$ are mutually independent and independent
of the problem instance.  We let $\Filt_t$ be the $\sigma$-algebra generated
by the initialization and the directions $u_i(\tau)$ for all $\tau<t$.
Consequently, $x_i(t)$ is $\Filt_t$-measurable and $u_i(t)$ is independent
of $\Filt_t$. 
 
With this notation, the update rule of Algorithm~\ref{alg:ozodgd}
admits the following compact stacked form
\begin{equation}\label{eq:stacked}
  \bx(t+1) \;=\; \bP\,\bx(t) \;-\; \eta_t\,\bg(t).
\end{equation}

\section{Regret and Convergence Guarantees}
\label{sec:theory}
\subsection{Online Regret Analysis}
\label{sec:regret}

We first establish the performance of Algorithm~\ref{alg:ozodgd} under the online setting. Our analysis combines the smoothed-function technique of zeroth-order online optimization, embodied in Lemma~\ref{lem:ball}, with the consensus-error analysis of distributed online optimization.

\begin{theorem}\label{thm:regret}
Run Algorithm~\ref{alg:ozodgd} with a constant step-size $\eta_t=\eta>0$ and arbitrary smoothing radii $\{\beta_t\}$,
starting from a consensual initialization $x_i(0)=\xbar(0)$. Let $R>0$ satisfy $\norm{\xbar(0)-x^\star}\leq R$. Then, under
Assumptions~\ref{as:network} and~\ref{as:reg}, for every
$j\in[N]$,
\begin{align}
 \E[\Reg_T(j)]
 \le{} &
N\eta TL_0^2
\left(\frac{N\!+\!d\!+\!4}{2N}+
\frac{3(\ln N\!+\!3)}{1\!-\!\sigma}+\frac{3\sqrt{d\!+\!4}}{\sqrt{1\!-\!\sigma^2}}\right) \notag\\
& + \frac{NR^2}{2\eta}
 +2NL_0\sum_{t=0}^{T-1}\beta_t.\label{eq:thm-tv}
\end{align}
\end{theorem}
\begin{proof}
See Appendix~\ref{app:regret}.
\end{proof}

\begin{corollary}\label{cor:const}
Under the conditions of Theorem~\ref{thm:regret}, the choice
\begin{equation*}
  \eta^2=\frac{R^2}{L_0^2T}\cdot\min\!\left\{
  \frac{N}{d},\frac{1-\sigma}{\ln N+1},\frac{\sqrt{1-\sigma}}{\sqrt{d}}
  \right\},
  \quad
  \beta_t\le\frac{R}{\sqrt T}
\end{equation*}
yields
\begin{equation}\label{eq:order}
  \E\bigl[\Reg_T(j)\bigr]
  =O\!\left(
  N\sqrt{T}\left[\max\!\left\{
\frac{d}{N},\frac{\ln N\!+\!1}{1\!-\!\sigma},\frac{\sqrt{d}}{\sqrt{1\!-\!\sigma}}
\right\}\right]^{\!\frac{1}{2}}
  \right)\!.
\end{equation}
\end{corollary}

Corollary~\ref{cor:const} establishes an $O(\sqrt{T})$ regret bound for the online setting with convex losses. While our analysis roughly follows standard existing approaches, we point out that the regret bound~\eqref{eq:order} exhibits improved dependences on $d$ and $N$ compared to the related existing works~\cite{NEURIPS2020_b3d6e130,NEURIPS2022_dececdcb,yang2026decentralized} in the corresponding static-network, uncompressed settings, under the respective assumptions. We note that this improvement is due to refined analysis rather than the use of ball-uniform estimators.

\subsection{Offline Convergence under Strong Convexity and Smoothness}
\label{sec:offline-convergence}

We next study the offline setting \eqref{eq:offline-model}, in which
$f_{i,t}\equiv f_i$. Here we consider the situation where each $f_i$ is strongly convex. Since a differentiable strongly convex function on
$\R^d$ cannot be globally Lipschitz, we replace Assumption~\ref{as:reg} by the following assumption throughout this subsection.

\begin{assumption}
\label{as:offline-sc-smooth}
For every $i\in[N]$, the fixed local objective $f_i:\R^d\to\R$ is
continuously differentiable, $\mu$-strongly convex, and $L_1$-smooth on
$\R^d$: For every $x,y\in\R^d$,
\[
\begin{gathered}
 f_i(y)\ge f_i(x)+\inner{\grad f_i(x)}{y-x}
                 +\frac{\mu}{2}\norm{y-x}^2,\\
 \norm{\grad f_i(x)-\grad f_i(y)}\le L_1\norm{x-y}.
\end{gathered}
\]
\end{assumption}

Under Assumption~\ref{as:offline-sc-smooth}, $f$ has a unique minimizer
$x^\star$. The problem is unconstrained, and thus $\grad f(x^\star)=0$. We also have the standard relation
$\mu\le L_1$.
Then, we let $f_i^\star$ be the infimum of $f_i$ over $\R^d$, which is guaranteed to be finite by the strong convexity of $f_i$. Define
\begin{equation}\label{eq:offline-heterogeneity}
  \mathcal H_f
  \coloneq \frac1N\sum_{i=1}^N\bigl(f_i(x^\star)-f_i^\star\bigr)\ge0.
\end{equation}
We further define the Lyapunov quantity
\begin{equation}\label{eq:offline-V}
  V_t
  \coloneq \norm{\xbar(t)-x^\star}^2+\frac{24L_1^2\eta_t}{\mu(1-\sigma^2)N}\norm{\xperp(t)}^2.
\end{equation}

\begin{theorem}
\label{thm:offline-convergence}
Suppose Assumptions~\ref{as:network} and
\ref{as:offline-sc-smooth} hold. Let
\[
  \eta_t=\frac{8}{\mu(t+b)},
  \qquad
  0<\beta_t\le\frac{\beta_{\mathrm{off}}}{t+b},
\]
where $\beta_{\mathrm{off}}>0$ and $b\ge \frac{256(d+1)L_1^2}{\mu^2(1-\sigma^2)}$.
Then, for all $0\le t\le T$,
\begin{equation}\label{eq:offline-V-rate}
  \E\!\left[\norm{\xbar(t)-x^\star}^2\right]\le \E[V_t]\le\frac{C_1}{t+b} + \frac{C_2}{(t+b)^2} + V_0\!\left(\frac{b}{t+b}\right)^{\!4},
\end{equation}
where
\begin{equation*}
\begin{aligned}
 C_{1}= \frac{512(d+1)L_1\mathcal{H}_f}{3\mu^2 N},
\end{aligned}
\end{equation*}
and $C_2$ is a constant of the order $O\!\left(
\frac{dL_1^3\mathcal{H}_f}{\mu^4(1-\sigma^2)^2}+\beta_{\mathrm{off}}^2L_1^2/\mu^2
\right)$. Furthermore, assuming consensual initialization $x_i(0)=\xbar(0)$, we have
\begin{align*}
& \frac1N\sum_{i=1}^N
 \E\norm{x_i(t)-\xbar(t)}^2 \le \frac{C_{\mathrm{con}}}{(t+b)^2},
\end{align*}
where $C_{\mathrm{con}}$ is a constant of the order $O\!\left(\frac{dL_1(\mathcal{H}_f+L_1V_0)}{\mu^2(1-\sigma^2)^2}\!+\!\beta_{\mathrm{off}}^2\right)$.
\end{theorem}

The detailed proof of Theorem~\ref{thm:offline-convergence} is postponed to Appendix~\ref{app:proof-offline-convergence}. We note that the rates established in Theorem~\ref{thm:offline-convergence} are consistent with the existing works~\cite{tang2021distributed,yi2022zeroth}. Especially, the $O(1/t)$ rate and linear dimension dependence in Theorem~\ref{thm:offline-convergence} are consistent with \cite{tang2021distributed,yi2022zeroth}, while the leading $O(d/(Nt))$ dependence matches the linear-speedup result in~\cite{yi2022zeroth}. For fixed $L_1,\mu$, and $\mathcal H_f$, the coefficient $C_1$ scales as $(d+1)/N$ and is independent of the network parameter $\sigma$.

\section{Differential Privacy Guarantees}
\label{sec:privacy}
\subsection{Single-Round Privacy}
We first isolate one agent in one round and derive the density of
the ball-uniform directional-derivative estimator
(Lemma~\ref{lem:density}), reduce the likelihood ratio between the
gradients of adjacent problems to a one-dimensional Cauchy variable
(Lemma~\ref{lem:cauchy}), and bound the mass of the high-risk region
(Theorem~\ref{thm:single}). The analytic proofs are collected in
Appendix~\ref{app:privacy}.

Following the discussion on the directional-derivative estimator in Section~\ref{sec:prelim_zo}, we let
$\bx^{\mathrm{lim}}(0)=\bx(0)$ and write
$\bx^{\mathrm{lim}}(t)\coloneq \mathrm{col}(x_i^{\mathrm{lim}}(t))_{i\in[N]}$.
Define
\[
  Y_i^{\mathrm{lim}}(t)\coloneq 
  \bigl(u_i(t)^{\!\top}\grad f_{i,t}(x_i^{\mathrm{lim}}(t))\bigr)u_i(t),
  \quad
  \tilde\eta_t\coloneq (d+2)\eta_t,
\]
and let
$\mathbf Y^{\mathrm{lim}}(t):=\mathrm{col}(Y_i^{\mathrm{lim}}(t))_{i\in[N]}$.
The limit-form recursion of Algorithm~\ref{alg:ozodgd} with the directional-derivative estimator is then
\begin{equation}\label{eq:limit-update}
  \bx^{\mathrm{lim}}(t+1)
  =\bP\bx^{\mathrm{lim}}(t)-\tilde\eta_t\mathbf Y^{\mathrm{lim}}(t).
\end{equation}
The adversary's observed trajectory is
\[
  \Pi_T^{\mathrm{lim}}
  \coloneq \bigl(x_i^{\mathrm{lim}}(t)\bigr)_{i\in[N],\,t\in\mathcal T_+}.
\]
The deterministic factor $d+2$ has been absorbed into the step size and
therefore does not affect likelihood ratios. All privacy statements below
refer to \eqref{eq:limit-update} and $\Pi_T^{\mathrm{lim}}$.

The starting point is the exact density of the estimator.

\begin{lemma}\label{lem:density}
Let $g\in\R^d\setminus\{0\}$ be deterministic, and $u\sim\Unif(\Ball^d)$. Then the random vector $(u^{\!\top}g)u$ is absolutely continuous with density
\begin{equation}\label{eq:density}
 \phi_g(y)=\begin{cases}
 \dfrac{(g^{\!\top}y)^{-d/2}}{\vol(\Ball^d)},&y\in S_g\setminus\{0\},\\[4pt]
 0,&y\notin S_g\ \text{or}\ y=0,
 \end{cases}
\end{equation}
where $S_g\coloneq \{y:g^{\!\top}y\ge\norm{y}^2\}$. 
\end{lemma}

The support $S_g$ of the probability law $\phi_g(y)$ is the solid ball with center $g/2$ and radius
$\norm{g}/2$. For $g,g'\in\R^d\setminus\{0\}$, the likelihood ratio on the
common positive-density region satisfies
\begin{equation}\label{eq:llr}
\begin{gathered}
 \frac{\phi_g(y)}{\phi_{g'}(y)}
 =\left(\frac{g'^{\!\top}y}{g^{\!\top}y}\right)^{\!d/2}\quad 
 \text{for a.e. }y\in(S_g\cap S_{g'})\setminus\{0\}.
\end{gathered}
\end{equation}
The privacy loss is therefore governed by the scalar ratio
$g'^{\top}y/g^{\!\top}y$. For any pair $(g,g')\in \R^d\times\R^d$, define the
\emph{relative gradient discrepancy}
\begin{equation}\label{eq:relative-discrepancy}
 \Dt(g,g')
 \coloneq 
 \begin{cases}
  0, & g=g', \\[2pt]
  \dfrac{2\norm{g-g'}}{\norm{g+g'}},&
        g+g'\ne0\ \text{and}\ g\ne g',\\[5pt]
  +\infty,& g+g'=0\ \text{and}\ g\ne g'.
 \end{cases}
\end{equation}
For a nonidentical pair with finite discrepancy, the next lemma shows that 
the single-step analysis can be reduced to one Cauchy variable.

\begin{lemma}\label{lem:cauchy}
Let $d\ge2$, $\epsilon>0$, and
$g,g'\in\R^d\setminus\{0\}$ satisfy $g\ne g'$ and $g+g'\ne0$.
For $u\sim\Unif(\Ball^d)$, denote $\varrho =\norm{u}$ and, on the probability-one event $\{u\ne0\}$,
$v = u/\norm{u}$. On the probability-one event
$\{(g+g')^{\!\top}v\ne0,\ g^{\!\top}v\ne0\}$, define
\[
  W = \frac{2(g'-g)^{\!\top}v}{(g+g')^{\!\top}v},\qquad
  X = \frac{g'^{\!\top}v}{g^{\!\top}v},
\]
and assign arbitrary values to $W$ and $X$ on the null complement. Then:
\begin{enumerate}[label=\textup{(\roman*)},leftmargin=2.4em,itemsep=2pt]
  \item Let $\mut=\Dt(g,g')\cos\theta$ and
        $\gat=\Dt(g,g')|\sin\theta|$, where $\theta$ is the angle between
        $g'-g$ and $(g+g')/2$. If $\gat>0$, then
        $W\sim\mathrm{Cauchy}(\mut,\gat)$; if $\gat=0$, then
        $W=\mut$ almost surely. In either case,
        $W$ is independent of $\varrho$.
  \item Almost surely, $X=(2+W)/(2-W)$, and
\[
\begin{aligned}
 \{X>e^{2\epsilon/d}\}&=\{2\tanh(\epsilon/d)<W<2\},\\
 \{X\le0\}&=\{|W|\ge2\},\\
 \{0<X<1\}&=\{-2<W<0\},
\end{aligned}
\]
        where the event equalities are understood modulo null events.
\end{enumerate}
\end{lemma}
 
The three events in Lemma~\ref{lem:cauchy}(ii) are exactly the three ways
a broadcast can betray the perturbed objective: The likelihood ratio~\eqref{eq:llr} can exceed $e^{\epsilon}$ inside the shared support ($X>e^{2\epsilon/d}$), or the
estimate can fall outside the support of the adjacent density, either
directionally ($X\le0$) or radially ($\varrho^2>X$). Because $W$ is a Cauchy variable
(or a degenerate point mass) with parameters satisfying
$\mut^2+\gat^2=\Dt(g,g')^2$, each event admits a closed-form
tail estimate depending only on $\Dt(g,g')$; carrying this out
(Appendix~\ref{app:privacy}) yields the single-step guarantee.
 
\begin{theorem}\label{thm:single}
Let $d\ge2$, $\epsilon>0$, and $g,g'\in\R^d$. For $h\in\R^d$, let
$\mathsf P_h$ denote the law of $(u^\top h)u$ for
$u\sim\Unif(\Ball^d)$; in particular, $\mathsf P_0$ is the point mass at
zero. Suppose $\Dt(g,g')<2\tanh(\epsilon/d)$. Then, for every measurable
$A\subseteq\R^d$,
\begin{equation}\label{eq:single-round-dp}
\begin{aligned}
 \mathsf P_g(A)
 &\le e^\epsilon\mathsf P_{g'}(A)
       +\delta_{\mathrm{step}}(\epsilon;\Dt(g,g')),\\
 \mathsf P_{g'}(A)
 &\le e^\epsilon\mathsf P_g(A)
       +\delta_{\mathrm{step}}(\epsilon;\Dt(g,g')),
\end{aligned}
\end{equation}
where $\delta_{\mathrm{step}}(\epsilon;0)\coloneq 0$ and, for $\chi\in(0,2\tanh(\epsilon/d))$,
\begin{equation}\label{eq:deltastep}
\begin{aligned}
 \delta_{\mathrm{step}}(\epsilon;\chi)
 \coloneq{} & {\frac1\pi}\arcsin\frac{\chi}{2\tanh(\epsilon/d)}
 +\frac{4\chi/\pi}{4-\chi^2}\\
 &+\frac d2\chi\left[1+\frac1\pi\ln\frac{2+\chi}{\chi}\right].
\end{aligned}
\end{equation}
Moreover, with $\epsilon$ fixed, $\delta_{\mathrm{step}}(\epsilon;\chi)$ is increasing over
$\chi\in[0,2\tanh(\epsilon/d))$ and tends to zero as $\chi\downarrow0$.
\end{theorem}

Interpreting the vectors $g,g'$ as the gradients of the distinct functions $f,f'$ in the adjacent problem instances, we can see from the condition $\Dt(g,g')<2\tanh(\epsilon/d)$ and the conclusion~\eqref{eq:single-round-dp} of Theorem~\ref{thm:single} that, to ensure single-round $(\epsilon,\delta)$-differential privacy via Theorem~\ref{thm:single}, the following assumption should be imposed.

\begin{assumption}\label{as:rel}
There exists $\Dt_{\max}\in[0,2)$ such that, for every instance pair adjacent
under Definition~\ref{def:adj} or Definition~\ref{def:offline-adj}, the
two distinct local objectives $f,f'$ at the affected agent (and at the affected round under the online setting) satisfy
\begin{equation}\label{eq:relative-gradient-assumption}
  \Dt\bigl(\grad f(x),\grad f'(x)\bigr)\le\Dt_{\max}
  \qquad\text{for all }x\in\R^d.
\end{equation}
\end{assumption}

In the following subsections, we will derive differential privacy guarantees for the whole observed trajectory $\Pi_T^{\mathrm{lim}}$ under Assumption~\ref{as:rel}.

\begin{remark}
\label{rem:intrinsic-privacy-scope}
Gupta et al.~\cite{gupta2025inherent} showed that the zero-preserving
zeroth-order oracles studied therein cannot provide a uniform privacy guarantee. Their hard instances compare an
identically zero objective with a nonzero quadratic or linear objective, so one estimator law contains a point mass at zero while the other is atomless.
Ball-uniform sampling alone does not remove this obstruction, since the
present limit estimator is also zero preserving. The distinction between  \cite{gupta2025inherent} and our result lies in
Assumption~\ref{as:rel}: For every \(g'\ne0\), we have $\Dt(0,g')=2$, and
hence the hard instances (i.e., the zero-nonzero pairs) proposed by \cite{gupta2025inherent} are outside
the situations considered by Assumption~\ref{as:rel} here.
Accordingly, our result should be interpreted as an (approximate) differential privacy guarantee for a restricted class of relatively close adjacent gradients, for which ball-uniform sampling yields comparable full-dimensional estimator densities (when the gradients are nonzero), rather than a universal inherent-privacy guarantee.

We also note that the condition~\eqref{eq:relative-gradient-assumption} necessarily requires the functions $f,f'$ to have the same stationary-point set and, under convexity, the same minimizer set, which might be restrictive. One possible approach to remedy this issue is to relax Assumption~\ref{as:rel} so that \eqref{eq:relative-gradient-assumption} is only required to hold on a high-probability region. For example, in the online setting, we may replace~\eqref{eq:relative-gradient-assumption} by
\[
\Dt\bigl(\grad f(x),\grad f'(x)\bigr)\le\Dt_{\max}
  \qquad\text{for}\ \  x\in K,
\]
where $K=K(F,F')\subseteq\R^d$ satisfies $\Prob(x_{i_0}^{\mathrm{lim}}(t_0)\notin K)\le q$ for some small $q>0$ ($i_0$ and $t_0$ denote the affected agent and round under adjacency). In this article, we focus on the simpler Assumption~\ref{as:rel}, and leave the relaxation and associated analysis for future work.
\end{remark}

\subsection{Privacy for the Online Setting}

\begin{revision}
\begin{theorem}\label{thm:privacy}
Let $d\ge2$ and $\epsilon>0$, and suppose Assumption~\ref{as:rel} holds
for the online adjacent problem instances with
$\Dt_{\max}<2\tanh(\epsilon/d)$. Then, for every horizon $T$, $\Pi_T^{\mathrm{lim}}$ is
$(\epsilon,\delta)$-differentially private, where
\begin{equation}\label{eq:main-delta}
  \delta
  =
  \min\left\{1,\,
  \delta_{\mathrm{step}}\bigl(\epsilon;\Dt_{\max}\bigr)\right\}.
\end{equation}
In particular, the bound is independent of $T$ and $N$.
\end{theorem}

\begin{proof}
If $\delta_{\mathrm{step}}(\epsilon;\Dt_{\max})\ge1$, the claim with
$\delta=1$ is immediate. Otherwise, fix adjacent online instances whose
unique differing function is indexed by $(i_0,t_0)$. Fix a round $t$ and a
deterministic candidate history of broadcast quantities $h_t=(h(0),\ldots,h(t))$ with $h(s)\in\R^{Nd}$ for $s\leq t$, and feed it to the
algorithm's iterations up to round $t$; let $h_i(t)$ be the $i$th block of $h(t)$.

If $t\ne t_0$, both instances use the same round-$t$ local objectives at every
agent, so their estimators' laws at round $t$ agree. If $t=t_0$, Assumption~\ref{as:rel} gives
\[
 \Dt\!\left(
   \grad f_{i_0,t_0}(h_{i_0}(t_0)),
   \grad f'_{i_0,t_0}(h_{i_0}(t_0))
 \right)\le\Dt_{\max}.
\]
Therefore, Theorem~\ref{thm:single} and monotonicity of
$\delta_{\mathrm{step}}(\epsilon;\cdot)$ imply that the two laws of
$Y_{i_0}^{\mathrm{lim}}(t_0)$ are
$(\epsilon,\delta_{\mathrm{step}}(\epsilon;\Dt_{\max}))$-indistinguishable.
The fresh directions are mutually independent and independent of the fixed
history; the estimators of all unaffected agents have the same laws under both
instances. Adjoining those estimators preserves the privacy parameters, and
\[
  \bP h(t)-\tilde\eta_t\mathbf Y^{\mathrm{lim}}(t)
\]
is deterministic post-processing. Thus all pointwise conditional round
mechanisms are $(0,0)$-private except round $t_0$, which has the displayed
parameters. Adaptive sequential composition~\cite{dwork2014algorithmic}, followed by appending the public,
data-independent initial state, proves the result for
$\Pi_T^{\mathrm{lim}}$.
\end{proof}
\end{revision}

\subsection{Privacy for the Offline Setting via Advanced Composition}
\label{sec:offline-privacy}

Under the offline setting \eqref{eq:offline-model}, replacing one
local objective changes the estimators' conditional distributions at every round.

\begin{lemma}[Basic and advanced adaptive composition]
\label{lem:adaptive-ac}
Let $Z_0,\ldots,Z_{T-1}$ be generated sequentially as follows: At round~$t$,
after a preceding output history $z_{<t}$ has been fixed, the algorithm
applies a randomized mechanism $\mathcal M_t(\,\cdot\,;z_{<t})$.
Suppose that, for every $t$, every fixed history $z_{<t}$, and every pair of
adjacent inputs, the map $F\mapsto\mathcal M_t(F;z_{<t})$ is
$(\epsilon_0,\delta_0)$-differentially private, where
$\epsilon_0\ge0$ and $\delta_0\in[0,1]$. Then basic composition gives
\[
  (Z_0,\ldots,Z_{T-1})
  \text{ is }
  \bigl(T\epsilon_0,\min\{1,T\delta_0\}\bigr)
  \text{-differentially private}.
\]
Moreover, for every $\widehat\delta\in(0,1)$, it is
$(\epsilon_T,\delta_T)$-differentially private, where
\begin{align}
  \epsilon_T
  &=\min\!\left\{
       T\epsilon_0,\,
       \epsilon_0\sqrt{2T\ln(1/\widehat\delta)}
       +T\epsilon_0\bigl(e^{\epsilon_0}-1\bigr)
     \right\},
  \label{eq:ac-epsilon}\\
  \delta_T
  &=\min\{1,T\delta_0+\widehat\delta\}.
  \label{eq:ac-delta}
\end{align}
\end{lemma}

Lemma~\ref{lem:adaptive-ac} follows from basic and advanced adaptive
composition; see \cite[Theorems~3.16 and~3.20]{dwork2014algorithmic}.
The minimum in \eqref{eq:ac-epsilon} is valid because the additive parameter $\delta$ in the $(\epsilon,\delta)$-differential privacy obtained by basic adaptive composition can be relaxed to $\min\{1,T\delta_0+\widehat\delta\}$, after which the smaller of the two available $\epsilon$ bounds from basic and advanced adaptive composition may be selected.

Applying Lemma~\ref{lem:adaptive-ac} to the offline setting gives its differential privacy guarantee, summarized in the following theorem.

\begin{theorem}
\label{thm:offline-advanced}
Let $d\ge2$ and fix $\epsilon_0>0$. Suppose
Assumption~\ref{as:rel} holds for the offline adjacent problem instances, with
\begin{equation*}
  \Dt_{\max}<2\tanh(\epsilon_0/d).
\end{equation*}
Define
\begin{equation}\label{eq:offline-delta0}
  \delta_0:=\min\!\left\{1,
  \delta_{\mathrm{step}}(\epsilon_0;\Dt_{\max})\right\}.
\end{equation}
Then, for every $\widehat\delta\in(0,1)$, the complete broadcast trajectory
$\Pi_T^{\mathrm{lim}}$ is
$(\epsilon_{\mathrm{off}},\delta_{\mathrm{off}})$-differentially private, where
\begin{align}
  \epsilon_{\mathrm{off}}
  &=\min\!\left\{
       T\epsilon_0,\,
       \epsilon_0\sqrt{2T\ln(1/\widehat\delta)}
       +T\epsilon_0\bigl(e^{\epsilon_0}-1\bigr)
     \right\},
  \label{eq:offline-epsilon}\\
  \delta_{\mathrm{off}}
  &=\min\{1,T\delta_0+\widehat\delta\}.
  \label{eq:offline-delta}
\end{align}
\end{theorem}

\begin{proof}
If $T\delta_0+\widehat\delta\ge1$, then
$\delta_{\mathrm{off}}=1$ and the stated claim is immediate. Hence suppose
$T\delta_0+\widehat\delta<1$. In particular, $\delta_0<1$, and therefore
\eqref{eq:offline-delta0} equals the untruncated single-round bound.

Fix a round $t$ and a deterministic candidate history
$h_t=(h(0),\ldots,h(t))$, with $h(s)\in\R^{Nd}$, and feed this history to
the algorithm's iterations up to $t$. 
Let $h_i(t)$ denote the $i$th $d$-dimensional block of $h(t)$. For simplicity, we let $\Prob_F^{h_t}$ denote the probability over the
fresh round-$t$ directions after the previous iterates of all agents have been replaced by the fixed blocks of $h_t$; $\Prob_{F'}^{h_t}$ has the analogous meaning for the adjacent instance $F'$. For the
affected agent $i_0$, Assumption~\ref{as:rel} gives
\[
  \Dt\!\left(\grad f_{i_0}(h_{i_0}(t)),
              \grad f'_{i_0}(h_{i_0}(t))\right)
  \le\Dt_{\max}.
\]
By Theorem~\ref{thm:single} and the monotonicity of
$\delta_{\mathrm{step}}(\epsilon_0;\cdot)$, for every measurable
$A\subseteq\R^d$,
\begin{equation}\label{eq:offline-one-round-conditional}
 \Prob_F^{h_t}\!\left(Y_{i_0}^{\mathrm{lim}}(t)\in A\right)
 \le e^{\epsilon_0}
 \Prob_{F'}^{h_t}\!\left(Y_{i_0}^{\mathrm{lim}}(t)\in A\right)
 +\delta_0,
\end{equation}
and the same inequality holds with $F$ and $F'$ interchanged.

For every $i\ne i_0$, the adjacent problems have the same local objectives,
so the corresponding estimators' laws agree. The random directions produced at round $t$ are mutually independent and independent of the fixed history. Hence adjoining the estimators of the unaffected agents preserves the privacy parameters. Thus,
for every measurable $B\subseteq\R^{Nd}$,
\[
 \Prob_F^{h_t}\!\left(\mathbf Y^{\mathrm{lim}}(t)\in B\right)
 \le e^{\epsilon_0}
 \Prob_{F'}^{h_t}\!\left(\mathbf Y^{\mathrm{lim}}(t)\in B\right)
 +\delta_0.
\]
The reverse inequality follows in the same way.
Finally, applying deterministic post-processing by the affine map $y\mapsto\bP h(t)-\tilde\eta_t y$ to the above inequality shows that, for every measurable
$C\subseteq\R^{Nd}$,
\[
\begin{aligned}
 &\Prob_F^{h_t}\!\left(
   \bP h(t)-\tilde\eta_t\mathbf Y^{\mathrm{lim}}(t)\in C\right)\\
 \le{} & e^{\epsilon_0}\Prob_{F'}^{h_t}\!\left(
   \bP h(t)-\tilde\eta_t\mathbf Y^{\mathrm{lim}}(t)\in C\right)
   +\delta_0.
\end{aligned}
\]
Thus the hypothesis of Lemma~\ref{lem:adaptive-ac} holds at every round, with $\mathcal{M}_t$ representing the Algorithm's iteration at round $t$. Applying the lemma to the complete trajectory of length $T$ and appending the public, data-independent initial state proves the claimed guarantees for $\Pi_T^{\mathrm{lim}}$.
\end{proof}

\section{Numerical Experiments}
\label{sec:experiments}

In this section, we first provide additional illustrations of the support geometry that motivates the design in
Section~\ref{sec:design}. Then we conduct a reconstruction-attack experiment, comparing the sphere-uniform estimator and the ball-uniform estimator when an adversary that observes the broadcast trajectory \eqref{eq:obs} and attempts to infer the private information is present.

\subsection{Support Geometry in Three Dimensions}

\begin{figure}[t]
  \centering
  \includegraphics[width=\columnwidth]{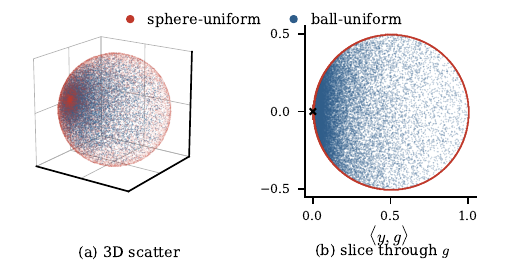}
  \caption{Normalized directional-derivative estimates $(u^{\!\top}g)u$ in $d=3$
  for a unit gradient $g$. (a)~Three-dimensional scatter of $2\times10^4$
  draws per law. (b)~Samples near a plane through the origin containing $g$.
  Sphere-uniform estimates lie on the boundary of $S_g$;
  ball-uniform estimates fill $S_g$.}
  \label{fig:support3d}
\end{figure}

Figure~\ref{fig:support3d} complements Figure~\ref{fig:supports} in $d=3$.
For a fixed unit gradient $g$, we draw directions from $\Unif(\Sphere^{2})$
and $\Unif(\Ball^{3})$ and form the directional-derivative estimators $d(u^{\!\top}g)u$ and
$(d+2)(u^{\!\top}g)u$, respectively. To compare the two laws on a common
scale, each output is normalized by its calibration constant
$c_{\mathcal D}$, which leaves $(u^{\!\top}g)u$ as in
Lemma~\ref{lem:density}. After this normalization, every sphere-uniform
estimate lies on the boundary sphere
$\norm{y-\tfrac12 g}=\tfrac12\norm{g}$ of the support $S_g$ in
\eqref{eq:density}, whereas the ball-uniform estimates fill $S_g$ and
concentrate toward the origin, in agreement with the density
$\phi_g$. The sphere-uniform estimator law is thus confined to the boundary of the support of the
ball-uniform estimator law, a set of zero Lebesgue measure.

\subsection{Reconstruction from the Broadcast Trajectory}

We consider the offline setting \eqref{eq:offline-model} with $N=10$
agents on a ring, $P_{ii}=P_{i,i\pm1}=1/3$, and $x_i(0)=0$. Each agent holds
one private $8\times8$ image from the UCI handwritten digits
dataset~\cite{alpaydin1998digits}, flattened to $a_i\in[0,1]^{64}$, and
uses $f_i(x)=\tfrac12\norm{x-a_i}^2$, which satisfies
Assumption~\ref{as:offline-sc-smooth} with $\mu=L_1=1$ and $d=64$. The
algorithmic parameters follow the conditions of Theorem~\ref{thm:offline-convergence}. 
We compare three decentralized-gradient-descent-type methods: the first-order
(FO) method with $g_i(t)=\grad f_i(x_i(t))$, the zeroth-order (ZO) method with the sphere-uniform estimator~\eqref{eq:est-prelim}, and the ZO method with the ball-uniform estimator~\eqref{eq:est-agent}.

The adversary knows $P$, $\{\eta_t\}$, the sampling law, and a public prior
$\bar a$, the pixel-wise mean of $360$ images disjoint from the private
ones. From two consecutive broadcasts it recovers
\[
  h_i(t):=\frac{1}{\eta_t}\Bigl(\sum_{j}P_{ij}x_j(t)-x_i(t+1)\Bigr)=g_i(t).
\]
Under FO, $a_i=x_i(0)-h_i(0)$ is revealed after one round. Under a
zeroth-order update, $h_i(t)=c_{\mathcal D}\bigl(u^{\!\top}(x_i(t)-a_i)\bigr)u$
with $u=u_i(t)$, so the direction $v_t:=h_i(t)/\norm{h_i(t)}$ is revealed.
Writing $u=\varrho w$ with $\varrho=\norm{u}$ and $\norm{w}=1$, each round
gives the scalar equation in $a_i$:
\begin{equation}\label{eq:attack-eq}
\begin{aligned}
  v_t^{\!\top}a_i
     +\Bigl(1-\frac{c_{\mathcal D}\varrho^2}{d}\Bigr)
      \bigl|w^{\!\top}(x_i(t)-a_i)\bigr|
      =v_t^{\!\top}x_i(t)-\frac{\norm{h_i(t)}}{d}.
\end{aligned}
\end{equation}
Since $c_{\mathcal D}\E[\varrho^2]=d$ for both estimators, the second term on the left-hand side has
zero conditional mean. For the sphere-uniform estimator, however,
$\varrho\equiv1$ and $c_{\mathcal D}=d$, so the second term on the left-hand side of~\eqref{eq:attack-eq} vanishes identically:
Each round yields an exact linear equation in $a_i$, and $d$ rounds with linearly independent directions determine $a_i$. For the ball-uniform estimator, the radial factor $\varrho^2$ remains random and the associated equations remain noisy. The adversary stacks the first $k$ equations \eqref{eq:attack-eq}, with coefficient matrix $V_k=[v_0,\ldots,v_{k-1}]^{\!\top}\in\R^{k\times d}$, and solves them by least squares, taking the minimum-norm correction of $\bar a$ when $k<d$. Accuracy is measured by $\mathrm{NRMSE}:=\norm{\hat a_i-a_i}/\norm{a_i}$, averaged over agents and over $10$ independent runs, where $\hat a_i$ denotes the adversary's estimate of $a_i$ produced by applying least squares to the available equations~\eqref{eq:attack-eq}.

\begin{figure}[t]
  \centering
  \includegraphics[width=\columnwidth]{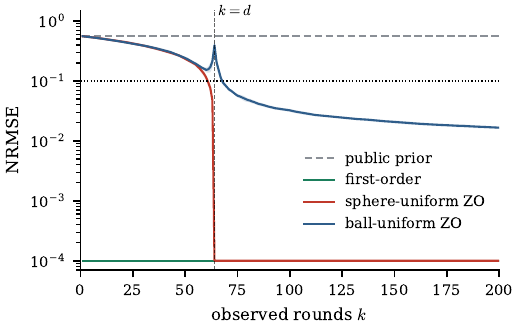}
  \caption{Reconstruction error of the private images from the broadcast
  trajectory versus the number of observed rounds ($N=10$, $d=64$; mean 
  over $10$ runs). Values below $10^{-4}$
  (first-order: about $10^{-16}$; sphere-uniform for $k\ge64$: about
  $10^{-12}$) are drawn at $10^{-4}$. Dotted line: $\mathrm{NRMSE}=0.1$;
  dashed line: $k=d$.}
  \label{fig:recon}
\end{figure}

\begin{figure}[t]
  \centering
  \includegraphics[width=\columnwidth]{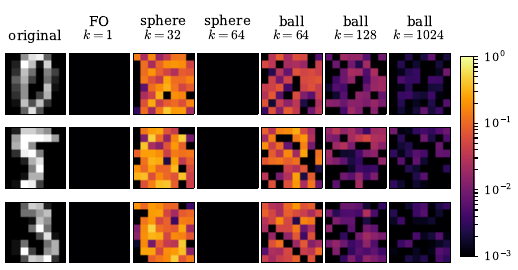}
  \caption{Private images of three agents (left) and the reconstruction errors $|\hat a_i-a_i|$ after $k$ observed rounds, on a
  logarithmic color scale. Errors below $10^{-3}$ are shown in black.}
  \label{fig:recon-images}
\end{figure}

Figure~\ref{fig:recon} reports the results. The FO update reveals every
$a_i$ after one round. For $k<d$, both zeroth-order estimators are limited
by the number of equations and behave similarly. At $k=d=64$, the sphere-uniform
reconstruction becomes exact up to rounding (NRMSE
$1.6\times10^{-12}$), whereas the ball-uniform error decreases only
gradually, to $2.4\times10^{-2}$ at $k=128$ and $5.9\times10^{-3}$ at
$k=1024$ (not shown in Figure~\ref{fig:recon}). Figure~\ref{fig:recon-images} shows the corresponding pixel-wise
errors for three agents: the sphere-uniform error vanishes at $k=64$,
whereas the ball-uniform error is still of order $10^{-1}$ there.

\begin{remark}
The transient peak of the ball-uniform curve at $k=d$ is a property of the
least-squares attack: It is essentially the
double-descent phenomenon of minimum-norm least
squares~\cite{belkin2019reconciling,hastie2022surprises}. Let
$e_i=a_i-\bar a$, let $\varepsilon=(\varepsilon_0,\ldots,\varepsilon_{k-1})$
collect the second terms on the left-hand side of~\eqref{eq:attack-eq}, and let $\Pi_k$ be the
orthogonal projector onto the row space of $V_k$. The reconstruction error
splits into two orthogonal parts,
\[
  \hat a_i-a_i=-(I-\Pi_k)e_i+V_k^{+}\varepsilon .
\]
The first part is the prior error in the directions not yet observed; its
expected squared norm is $(1-k/d)_+\norm{e_i}^2$, which vanishes at $k=d$.
The second part is the noise amplified by $V_k^{+}$; its energy loosely scales with
$\mathrm{tr}\bigl((V_k^{\top}V_k)^{+}\bigr)$, which for random directions is
approximately $dk/(d-k)$ for $k<d$ and $d^2/(k-d)$ for $k>d$. It therefore
can be very large near the threshold $k=d$, where $V_k$ becomes square
and can be poorly conditioned. The matrix associated with the sphere-uniform estimator is equally ill-conditioned at $k=d$, but its noise vanishes identically, so its error drops to zero instead of peaking.
\end{remark}

\section{Conclusion}
\label{sec:conclusion}

In this article, we studied distributed zeroth-order optimization with privacy leakage from iterative local variable exchange. We proposed exploiting the random query directions of zeroth-order estimation as a privacy source and combined a ball-uniform two-point gradient estimator with consensus updates. The algorithm achieves convergence rates comparable to those of existing methods in both online and offline settings. For privacy, we derive the single-round and trajectory-level \((\epsilon,\delta)\)-differential privacy guarantees. These results indicate that intrinsic zeroth-order randomness can serve as a promising alternative to additive-noise mechanisms.

Several directions remain open. First, it is promising to extend our framework to other distributed zeroth-order algorithms, in particular
gradient-tracking methods~\cite{qu2018harnessing,nedic2017diging,tang2021distributed},
whose auxiliary tracking variables are also exchanged and therefore enlarge
the adversary's observed trajectory. Second, the privacy guarantees established here concern the directional-derivative estimator; characterizing the distribution of the finite-difference estimator
would yield guarantees for the original algorithm itself. Third, constrained problems call for projected updates, and establishing privacy guarantees for the projected algorithm 
is a natural next step.

\paragraph*{Statement on AI use}
After finishing a first draft, the authors used ChatGPT 6 Astra to suggest possible refinements of the results and proofs. Notably, the AI system suggested using the Poincar\'{e} inequality on the sphere in proving Lemma~\ref{lem:ball}, and also pointed out several loose derivations. In addition, the AI system was used for polishing the language. All AI-assisted revisions were reviewed and verified by the authors, who take full responsibility for the manuscript's content and accuracy.

\appendices
\section{Proof of Lemma~\ref{lem:ball}}
\label{app:estimator}
The following lemma lists standard results in zeroth-order optimization; see, e.g., \cite{flaxman2005online,shamir2017optimal,malik2020derivative}.
\begin{lemma}\label{lem:smooth}
Let $h:\R^d\to\R$ be continuous. Then for every $s>0$, the
smoothed function $\hat h^{s}$ of \eqref{eq:smooth} is differentiable, and
\begin{equation}
\label{eq:smoothed_ver_gradient_sphere}
\begin{aligned}
\grad\hat h^{s}(x) ={} &
\frac{d}{2s}\,\E_{v}[(h(x+sv)-h(x-sv))v],
\end{aligned}
\end{equation}
where $v\sim\Unif(\Sphere^{d-1})$.
Furthermore,
\begin{enumerate}[label=\textup{(\roman*)},leftmargin=2.4em,itemsep=2pt]
\item If $h$ is convex, then $\hat h^{s}$ is convex.
\item If $h$ is $L_0$-Lipschitz, then $\hat h^{s}$ is also $L_0$-Lipschitz, and $\bigl|\hat h^{s}(x)-h(x)\bigr|\le s\,L_0,\,\forall x$.

\item If $h$ is $L_1$-smooth, then $\hat h^{s}$ is also $L_1$-smooth, and $\big\|\nabla\hat h^{s}(x)-\nabla h(x)\big\|\le s\,L_1,\,\forall x$.
\end{enumerate}
\end{lemma}

We now prove Lemma~\ref{lem:ball}.
Since $\Prob(u=0)=0$ when $u\sim\Unif(\Ball^d)$, for simplicity we shall just assume $u\neq 0$. Write $u=\varrho v$ with $\varrho=\norm{u}\in(0,1]$ and
$v=u/\norm{u}\in\Sphere^{d-1}$. For $u\sim\Unif(\Ball^d)$, it is known from probability theory that $\varrho$ and $v$ are independent, $v\sim\Unif(\Sphere^{d-1})$, and $\varrho$ has density
$p_{\varrho}(r)=d\,r^{d-1}$ over $r\in(0,1]$.

\medskip\noindent
\emph{(i)} We have
\[
G_h(x;\beta,u)
=\frac{d\!+\!2}{d}\varrho^2\cdot\frac{d}{2(\beta\varrho)}\bigl(h(x\!+\!\beta \varrho v)-h(x\!-\!\beta\varrho v)\bigr)v.
\]
Taking expectations over $v$ conditioned on $\varrho$ and invoking~\eqref{eq:smoothed_ver_gradient_sphere} at radius $s = \beta \varrho$, we get
\[
  \E\Econd*{G_h(x;\beta,u)}{\varrho}=\frac{d+2}{d}\,\varrho^2\,\grad\hat h^{\beta \varrho}(x).
\]
Averaging over $\varrho$ with density $d\,r^{d-1}$ gives
\begin{equation}
\label{eq:expectation_estimator_temp}
\begin{aligned}
 \E[G_h(x;\beta,u)]
 &=(d\!+\!2)\int_0^1r^{d+1}\grad\hat h^{\beta r}(x)\,\ud r\\
 &=\grad\!\left[(d\!+\!2)\int_0^1r^{d+1}\hat h^{\beta r}(x)\,\ud r\right]
=\grad h^\beta(x).
\end{aligned}
\end{equation}
Here the interchange of differentiation and integration can be justified by the dominated convergence theorem. 

\medskip\noindent
\emph{(ii)} When $h$ is convex, by Lemma~\ref{lem:smooth}(i), $\hat h^{\beta r}$ is convex. Since $r^{d+1}\ge 0$ for $r\in [0,1]$, $h^\beta$ is a nonnegative weighted integral of convex functions, hence convex.

\medskip\noindent
\emph{(iii)} When $h$ is $L_0$-Lipschitz, Lemma~\ref{lem:smooth}(ii) implies that $\hat{h}^{\beta r}$ is $L_0$-Lipschitz for all $r\in(0,1]$. As a result,
\[
\begin{aligned}
|h^\beta(x)\!-\!h^\beta(y)|
={} & \left|(d\!+\!2)\!\int_0^1 r^{d+1}(\hat{h}^{\beta r}(x)-\hat{h}^{\beta r}(y))\,\ud r\right| \\
\leq{} & (d\!+\!2)\!\int_0^1 r^{d+1}L_0\|x\!-\!y\|\,\ud r
\leq L_0 \|x\!-\!y\|,
\end{aligned}
\]
showing that $h^\beta$ is $L_0$-Lipschitz. Furthermore, using the bound $|\hat h^{\beta r}(x)-h(x)|\le\beta r L_0$ from Lemma~\ref{lem:smooth}(ii) and $(d+2)\int_0^1 r^{d+1}\,\ud r=1$, we get
\[
\begin{aligned}
 \bigl|h^\beta(x)-h(x)\bigr|
 ={} & \Bigl|(d+2)\int_0^1r^{d+1}\bigl(\hat h^{\beta r}(x)-h(x)\bigr)\,\ud r\Bigr|\\
 &\le(d+2)\int_0^1r^{d+1}\beta rL_0\,\ud r
\le\beta L_0.
\end{aligned}
\]
To show the bound on the second moment of $\norm{G_h(x;\beta,u)}$, denote for each $s>0$
\[
q_s(v) = \frac{h(x+sv)-h(x-sv)}{2s}.
\]
When $h$ is $L_0$-Lipschitz, it is straightforward to verify that $q_s\colon \R^d\to\R$ is also $L_0$-Lipschitz. As a result, the spherical gradient of $q_s$, which we temporarily denote by $\nabla_\Sphere q_s$, exists almost everywhere on $\Sphere^{d-1}$ and is integrable, satisfying $\norm{\nabla_\Sphere q_s(v)}\leq L_0$. Moreover, $\E_v[q_s(v)] = 0$ for $v\sim\Unif(\Sphere^{d-1})$. For $d\geq 2$, by applying the Poincar\'{e} inequality on the sphere (see, e.g., \cite[Theorem~1.1]{dolbeault2013sharp} and the discussion therein), we get
\[
\mathbb{E}_v\!\left[|q_s(v)|^2\right]
\le \frac{1}{d-1}\mathbb{E}_v\!\left[\norm{\nabla_\Sphere q_s(v)}^2\right] \leq\frac{L_0^2}{d-1},
\]
where $v\sim\Unif(\Sphere^{d-1})$. Now, recalling that $\varrho=\|u\|$ and $v=u/\|u\|$, by $G_h(x;\beta,u) = (d+2)\varrho^2\cdot q_{\beta\varrho}(v)\,v$, we obtain
\[
\begin{aligned}
\E\!\left[\norm{G_h(x;\beta,u)}^2\right]
={} & \E\!\left[\E\Econd*{(d\!+\!2)^2\varrho^4\cdot |q_{\beta\varrho}(v)|^2}{\varrho}\right] \\
\le{} & \E\!\left[
\mfrac{(d\!+\!2)^2L_0^2}{d\!-\!1}\varrho^4
\right] = \mfrac{d(d+2)^2}{(d\!-\!1)(d\!+\!4)}L_0^2 \\
\le{} & (d+4)L_0^2,
\end{aligned}
\]
where we used $\mathbb{E}[\varrho^4] = d/(d+4)$. When $d=1$, we have $|q_{\beta\varrho}(v)|\leq L_0$, which will lead to $\E\!\left[\norm{G_h(x;\beta,u)}^2\right]\leq 5L_0^2=(d+4)L_0^2$.

\medskip\noindent
\emph{(iv)} When $h$ is $L_1$-smooth, Lemma~\ref{lem:smooth}(iii) shows that $\hat h^{\beta r}$ is $L_1$-smooth for all $r\in(0,1]$. By~\eqref{eq:expectation_estimator_temp}, we get
\[
\begin{aligned}
& \norm{\nabla h^\beta(x)-\nabla h^\beta(y)} \\
={} &
\norm{(d+2)\int_0^1 r^{d+1}\left(\nabla \hat{h}^{\beta r}(x)-\nabla \hat{h}^{\beta r}(y)\right)\ud r} \\
\le{} &
(d+2)\int_0^1 r^{d+1}L_1\norm{x-y}\ud r \le L_1\|x-y\|,
\end{aligned}
\]
showing that $h^\beta$ is $L_1$-smooth. To bound $\norm{\nabla h^\beta(x)-\nabla h(x)}$, we use $\big\|\nabla\hat h^{\beta r}(x)-\nabla h(x)\big\|\le \beta r\,L_1$ from Lemma~\ref{lem:smooth}(iii) and obtain
\[
\begin{aligned}
& \norm{\nabla h^\beta(x)-\nabla h(x)} \\
={} & \norm{(d+2)\int_0^1 r^{d+1}\left(\nabla \hat{h}^{\beta r}(x)-\nabla {h}(x)\right)\ud r} \\
\le{} & (d+2)\int_0^1 r^{d+1} \beta r L_1\,\ud r\le \beta L_1.
\end{aligned}
\]
To bound the second moment of $\norm{G_h(x;\beta,u)}$, we denote $\tilde{q}_s(v) = q_s(v)-\inner{v}{\nabla h(x)}$ for each $s>0$. Since $h$ is assumed to be $L_1$-smooth, we have
\[
\begin{aligned}
\norm{\nabla_v\tilde{q}_s(v)}
\! ={} &
\!\norm{\frac{\nabla h(x\!+\!sv)\!-\!\nabla h(x)}{2} \!+\! \frac{\nabla h(x\!-\!sv)\!-\!\nabla h(x)}{2}} \\
\le{} & \!\frac{\norm{\nabla h(x\!+\!sv)\!-\!\nabla h(x)} \!+\! \norm{\nabla h(x\!-\!sv)\!-\!\nabla h(x)}}{2} \\
\le{} & \frac{L_1s\norm{v}+L_1s\norm{v}}{2} = L_1s\norm{v}.
\end{aligned}
\]
Thus $\norm{\nabla_\Sphere \tilde{q}_s(v)}\leq L_1s$ for $v\in\Sphere^{d-1}$. Moreover, it's evident that $\E_v[\tilde{q}_s(v)]=0$ when $v\sim\Unif(\Sphere^{d-1})$. For $d\geq 2$, we can then apply the Poincar\'{e} inequality on the sphere to obtain
\[
\E_v\left[|\tilde{q}_s(v)|^2\right]
\leq\frac{1}{d-1}\E_v\!\left[\norm{\nabla_\Sphere \tilde{q}_s(v)}^2\right] \leq \frac{L_1^2s^2}{d-1}
\]
with $v\sim\Unif(\Sphere^{d-1})$. Now for $u\sim\Unif(\Ball^d)$, we again decompose it as $u=\varrho v$ with $\varrho=\norm{u}$ and $v\sim\Unif(\Sphere^{d-1})$ being independent. We then have
\[
\begin{aligned}
& \E\left[\norm{G_h(x;\beta,u)-(d\!+\!2)\inner{u}{\nabla h(x)}u}^2\right] \\
={} & \E\left[\E\Econd*{\norm{(d\!+\!2)\varrho^2 \tilde{q}_{\beta\varrho}(v)\,v}^2}{\varrho}\right]
\le \E\!\left[(d\!+\!2)^2\varrho^4\cdot\mfrac{L_1^2\beta^2\varrho^2}{d\!-\!1}\right] \\
={} &
\mfrac{d(d+2)^2}{(d-1)(d+6)}\beta^2 L_1^2 \le (d+2)\beta^2 L_1^2.
\end{aligned}
\]
For $d=1$, it can be shown by a standard Taylor expansion argument that the above inequality also holds. Applying $\norm{a+b}^2\le 2\norm{a}^2+2\norm{b}^2$ then yields
\[
\begin{aligned}
& \E\!\left[\norm{G_h(x;\beta,u)}^2\right] \\
 \le{} & 
2\E\!\left[\norm{(d\!+\!2)\inner{u}{\nabla h(x)}u}^2\right]
+2(d\!+\!2)\beta^2L_1^2 \\
={} & 2(d\!+\!2)^2\E\!\left[\varrho^4 \nabla h(x)^\top vv^\top \nabla h(x)\right]
+2(d\!+\!2)\beta^2L_1^2 \\
={} &  2\mfrac{(d\!+\!2)^2}{d+4}\norm{\grad h(x)}^2
 +2(d\!+\!2)\beta^2 L_1^2,
\end{aligned}
\]
and using $(d+2)^2/(d+4)\leq d+1$ completes the proof.

\section{Proof of Theorem~\ref{thm:regret}}
\label{app:regret}

The following lemma is the standard online-gradient-descent bound stated for a generic descent sequence; see, e.g., \cite[Section 2.3]{nemirovski2009robust}.

\begin{lemma}\label{lem:ogd-tv}
Let $\eta>0$, and let
$\{z(t)\}_{t=0}^{T}\subset\R^d$ and $\{h(t)\}_{t=0}^{T-1}\subset\R^d$
satisfy $z(t+1)=z(t)-\eta\,h(t)$ for all $0\le t\le T-1$. Then
\[
  \sum\nolimits_{t=0}^{T-1}\inner{h(t)}{z(t) \!-\! x^\star}
  \le\mfrac{\|z(0)\!-\!x^\star\|^2}{2\,\eta}
        +\mfrac{\eta}{2}\sum\nolimits_{t=0}^{T-1}\norm{h(t)}^2.
\]
\end{lemma}

Our proof roughly follows the analysis of standard distributed online gradient descent, applied to the smoothed losses
$f^{\beta_t}_{i,t}$, i.e., the function $h^\beta$ defined by~\eqref{eq:Phi} associated with $h=f_{i,t}$ and $\beta=\beta_t$, with exact gradients replaced by gradient estimates $g_i(t)$. Lemma~\ref{lem:ball} guarantees that $\E\Econd{g_i(t)}{\Filt_t}=\nabla f^{\beta_t}_{i,t}(x_i(t))$. 
On the other hand, we leverage techniques from \cite[Appendix~A]{wan2025optimal} and conduct a refined analysis on the stochastic consensus error to obtain improved dependence of the regret on $N$ and $d$.

Since $P$ is doubly stochastic, the network average obeys
$\xbar(t+1)=\xbar(t)-\eta\gbar(t)$, while the consensus error obeys the recursion
$\xperp(t+1)=((P-J)\otimes I_d)\xperp(t)-\eta(M\otimes I_d)\bg(t)$.

\vspace{3pt}
\noindent{}\textbf{Step I: Bounding consensus error.} Denote $P^\perp_\tau = P^\tau-J$ for $\tau\geq 0$. Unrolling the consensus error recursion
from $\xperp(0)=0$ and noting $(P-J) M=P-J$ gives
\[
\xperp(t) = -\eta\sum\nolimits_{\tau=0}^{t-1}(P^\perp_{t-1-\tau} \otimes I_d) \bg(\tau).
\]
Denoting $\xi_i(\tau) = g_i(\tau)-\nabla f^{\beta_\tau}_{i,\tau}(x_i(\tau))$ for each $\tau$, we get
\begin{align}
x_j(t)-\xbar(t)
={} &
{-\eta}\sum\nolimits_{\tau=0}^{t-1}\sum\nolimits_{i}
(P^\perp_{t-1-\tau})_{j,i}\nabla f^{\beta_\tau}_{i,\tau}(x_i(\tau)) \notag\\
& -\eta\sum\nolimits_{\tau=0}^{t-1}\sum\nolimits_{i}
(P^\perp_{t-1-\tau})_{j,i}\xi_i(\tau),
\label{eq:consensus_error_bound1}
\end{align}
To bound the first term on the right-hand side of~\eqref{eq:consensus_error_bound1}, noting that each $f^{\beta_\tau}_{i,\tau}$ is $L_0$-Lipschitz by Lemma~\ref{lem:ball}, we have
\[
\begin{aligned}
& \left\|\sum\nolimits_{\tau=0}^{t-1}\sum\nolimits_{i}
(P^\perp_{t-1-\tau})_{j,i}\nabla f^{\beta_\tau}_{i,\tau}(x_i(\tau))\right\| \\
\leq{} & \sum\nolimits_{\tau=0}^{t-1}\sum\nolimits_{i} \left|(P^\perp_{t-1-\tau})_{j,i}\right| L_0
={} L_0 \sum\nolimits_{\tau=0}^{t-1}\left\|e_j^\top P^\perp_{\tau}\right\|_1 \\
\leq{} & L_0 \sum\nolimits_{\tau=0}^{t-1}\min\!\left\{2, \sqrt{N}\left\|e_j^\top P^\perp_{\tau}\right\|\right\} \\
\leq{} & L_0 \sum\nolimits_{\tau=0}^{t-1}\min\!\left\{2, \!\sqrt{N}\sigma^{\tau}\!\right\},
\end{aligned}
\]
where the third step follows because $P^{\tau}$ is doubly stochastic and hence
$\norm{e_j^\top P^\perp_{\tau}}_1
\leq{} \norm{e_j^\top P^{\tau}}_1 + N^{-1}\norm{\mathbf{1}}_1
\leq 2$,
while the fourth step follows by $\norm{e_j^\top P^\perp_{\tau}}\leq \norm{P^\tau-J}_2\leq \norm{P-J}_2^{\tau}\leq\sigma^{\tau}$. Furthermore, by letting $\check\tau = \max\{0,\lceil-\ln (\sqrt{N}/2)/\ln\sigma\rceil\}$ when $\sigma>0$, we get
\[
\begin{aligned}
& \sum\nolimits_{\tau=0}^{t-1}\!\min\!\left\{2, \!\sqrt{N}\sigma^{\tau}\!\right\}
\le \sum\nolimits_{\tau=0}^{\check\tau-1}2+\sum\nolimits_{\tau=\check\tau}^{\infty}\sqrt{N}\sigma^{\tau} \\
\leq{} & 2\max\Big\{0,\Big\lceil\mfrac{\ln(\sqrt{N}/2)}{-\ln\sigma}\Big\rceil\Big\}+\mfrac{\sqrt{N}\sigma^{\check\tau}}{1-\sigma}
\leq \mfrac{\ln N + 3}{1-\sigma},
\end{aligned}
\]
and the above bound obviously holds for $\sigma=0$. Therefore
\[
\norm{\sum\nolimits_{\tau=0}^{t-1}\sum\nolimits_{i}
(P^\perp_{t-1-\tau})_{j,i}\nabla f^{\beta_\tau}_{i,\tau}(x_i(\tau))}
\leq\mfrac{(\ln N+3)L_0}{1-\sigma}.
\]
To bound the second term on the right-hand side of~\eqref{eq:consensus_error_bound1}, we note that
\[
\begin{aligned}
& \norm{\sum\nolimits_{\tau=0}^{t-1}\sum\nolimits_{i}
(P^\perp_{t-1-\tau})_{j,i}\xi_i(\tau)}^2 \\
={} & \sum_{\tau,s=0}^{t-1}
\sum_{i,k}
(P^\perp_{t-1-\tau})_{j,i}(P^\perp_{t-1-s})_{j,k}
\langle\xi_i(\tau),\xi_k(s)\rangle.
\end{aligned}
\]
By taking the expectation on both sides and noting that
\[
\begin{aligned}
\E\!\left[\inner{\xi_i(\tau)}{\xi_k(s)}\right]
={} & \E\!\left[\E\Econd*{\inner{\xi_i(\tau)}{\xi_k(s)}}{\Filt_{\min\{\tau,s\}}}\right] \\
={} & \delta_{i,k}\delta_{\tau,s}
\E\big[\|\xi_i(\tau)\|^2\big],
\end{aligned}
\]
where $\delta_{i,k}$ denotes the Kronecker $\delta$ function, we get
\[
\begin{aligned}
& \E\norm{\sum\nolimits_{\tau=0}^{t-1}\sum\nolimits_{i}
(P^\perp_{t-1-\tau})_{j,i}\xi_i(\tau)}^2 \\
={} & \sum\nolimits_{\tau=0}^{t-1}\sum\nolimits_{i}
(P^\perp_{t-1-\tau})_{j,i}^2 \E\!\left[\|\xi_i(\tau)\|^2\right].
\end{aligned}
\]
Since $\xi_i(\tau)=g_i(\tau)-\E\Econd{g_i(\tau)}{\Filt_\tau}$, we have $\E\!\left[\|\xi_i(\tau)\|^2\right]\leq \E\!\left[\|g_i(\tau)\|^2\right]$, and by Lemma~\ref{lem:ball}, we further have $\E\!\left[\|g_i(\tau)\|^2\right]\leq L_0^2(d+4)$. Consequently,
\[
\begin{aligned}
& \E\norm{\sum\nolimits_{\tau=0}^{t-1}\sum\nolimits_{i}
(P^\perp_{t-1-\tau})_{j,i}\xi_i(\tau)}^2 \\
\leq{} & L_0^2(d\!+\!4)\sum_{\tau=0}^{t-1}\sum_{i}
(P^\perp_{t-1-\tau})_{j,i}^2
=
L_0^2(d\!+\!4)\sum_{\tau=0}^{t-1}\norm{e_j^\top P^\perp_{\tau}}^2 \\
\leq{} & L_0^2(d+4)\sum\nolimits_{\tau=0}^{t-1}\sigma^{2\tau}
\leq \mfrac{L_0^2(d+4)}{1-\sigma^2},
\end{aligned}
\]
and by Jensen's inequality,
\[
\E\norm{\sum\nolimits_{\tau=0}^{t-1}\sum\nolimits_{i}
(P^\perp_{t-1-\tau})_{j,i}\xi_i(\tau)}
\leq\mfrac{L_0\sqrt{d+4}}{\sqrt{1-\sigma^2}}.
\]
Summarizing the previous bounds, we get
\begin{equation}
\label{eq:consensus_error_bound}
\E\!\left[
\norm{x_j(t)-\xbar(t)}\right]
\leq \eta L_0\Big(\mfrac{\ln N+3}{1-\sigma}+\mfrac{\sqrt{d+4}}{\sqrt{1-\sigma^2}}\Big).
\end{equation}

\noindent{}\textbf{Step II: Bounding the regret.} We have
\begin{align}
 &f_{i,t}(x_j(t))-f_{i,t}(x^\star)
\le f^{\beta_t}_{i,t}(x_j(t)) - f^{\beta_t}_{i,t}(x^\star) + 2\beta_t L_0 \notag\\
\le{} & f^{\beta_t}_{i,t}(x_i(t)) \!-\! f^{\beta_t}_{i,t}(x^\star) + 2\beta_t L_0
+|f^{\beta_t}_{i,t}(x_j(t))\!-\!f^{\beta_t}_{i,t}(x_i(t))| \notag\\
 \le{} &\inner{\grad f^{\beta_t}_{i,t}(x_i(t))}{x_i(t)\!-\!x^\star}
       +L_0\norm{x_j(t)\!-\!x_i(t)} + 2\beta_tL_0 \notag\\
\le{} & \inner{\grad f^{\beta_t}_{i,t}(x_i(t))}{\xbar(t) \!-\!x^\star}
+\norm{\grad f^{\beta_t}_{i,t}(x_i(t))}\norm{x_i(t)\!-\!\xbar(t)} \notag\\
& + L_0\norm{x_j(t)\!-\!\xbar(t)} + L_0\norm{\xbar(t)\!-\!x_i(t)} + 2\beta_tL_0 \notag\\
\le{} & \inner{\grad f^{\beta_t}_{i,t}(x_i(t))}{\xbar(t) \!-\!x^\star}
+ L_0\norm{x_j(t)\!-\!\xbar(t)} \notag\\
& + 2L_0\norm{x_i(t)\!-\!\xbar(t)} + 2\beta_tL_0,
\label{eq:step2}
\end{align}
where we used Lemma~\ref{lem:ball}(iii) in the first step, used the convexity of $f^{\beta_t}_{i,t}$ in the third step, and used $L_0$-Lipschitz continuity of $f^{\beta_t}_{i,t}$ in the third and the last steps.
By Lemma~\ref{lem:ball}(i),
$\E\Econd*{g_i(t)}{\Filt_t}=\grad f^{\beta_t}_{i,t}(x_i(t))$. Summing \eqref{eq:step2}
over $i$ and $t$, taking expectations, and using~\eqref{eq:consensus_error_bound} gives
\[
\begin{aligned}
 \E[\Reg_T(j)]
 \le{} & N\E\!\left[\sum\nolimits_{t=0}^{T-1}\!\inner{\gbar(t)}{\xbar(t)\!-\!x^\star}\right]
 \!+2NL_0\!\sum\nolimits_{t=0}^{T-1}\!\beta_t\\
 &+3NL_0^2\eta T\bigg(\frac{\ln N+3}{1-\sigma}+\frac{\sqrt{d+4}}{\sqrt{1-\sigma^2}}\bigg)
\end{aligned}
\]
By applying Lemma~\ref{lem:ogd-tv} to the first term on the right-hand side, we get
\[
\E\!\left[\sum\nolimits_{t=0}^{T-1}\inner{\gbar(t)}{\xbar(t)-x^\star}\right]
\leq \frac{R^2}{2\eta}
+\frac{\eta}{2}\sum\nolimits_{t=0}^{T-1}\mathbb{E}\!\left[\norm{\gbar(t)}^2\right]\!,
\]
while Lemma~\ref{lem:ball}(iii) implies
\[
\begin{aligned}
& \mathbb{E}\Econd*{\norm{\gbar(t)}^2}{\Filt_t}
=\mfrac{1}{N^2}\mathbb{E}\Econd*{\norm{\sum\nolimits_{i}g_i(t)}^2}{\Filt_t} \\
={} & 
\mfrac{1}{N^2}\left(
\norm{\sum\nolimits_i\nabla f^{\beta_t}_{i,t}(x_i(t))}^2
+\sum\nolimits_i\E\Econd*{\norm{\xi_i(t)}^2}{\Filt_t}\right) \\
\le{} &
\mfrac{1}{N^2}\left(
N^2L_0^2 + N\cdot L_0^2(d+4)
\right) = L_0^2\left(1+N^{-1}(d+4)\right),
\end{aligned}
\]
where we recall that $\xi_i(t) = g_i(t) - \nabla f^{\beta_t}_{i,t}(x_i(t))$. Summarizing all the obtained bounds yields \eqref{eq:thm-tv}.

\section{Proof of Theorem~\ref{thm:offline-convergence}}
\label{app:proof-offline-convergence}

\noindent
\textbf{Step I: Second moment and bias.} Lemma~\ref{lem:ball}(iv) gives
\begin{equation}\label{eq:offline-estimator-moment}
\begin{aligned}
 \E\Econd*{\norm{g_i(t)}^2}{\Filt_t}
 \le{} & 2(d\!+\!1)\norm{\grad f_i(x_i(t))}^2
 +2(d\!+\!2)L_1^2\beta_t^2.
\end{aligned}
\end{equation}
Moreover, the $L_1$-smoothness of $f_i$ leads to
\[
\begin{aligned}
\norm{\grad f_i(x_i(t))}^2\le {} &
2\norm{\grad f_i(x_i(t))\!-\!\grad f_i(x^\star)}^2
+2\norm{\grad f_i(x^\star)}^2 \\
\le{} & 2L_1^2\norm{x_i(t)-x^\star}^2
+4L_1(f_i(x^\star)-f_i^\star),
\end{aligned}
\]
where we used $\norm{\grad f_i(x^\star)}^2\le2L_1(f_i(x^\star)-f_i^\star)$ (see, e.g., \cite[Section 1.2.3]{nesterov2018lectures}). Summing
\eqref{eq:offline-estimator-moment} over $i$ and using $\sum_i\norm{x_i(t)-x^\star}^2 = N\norm{\xbar(t)-x^\star}^2+\norm{\bx_\perp}^2$ yields
\begin{equation}\label{eq:offline-moment}
\begin{aligned}
 \sum\nolimits_{i}\!\E\Econd*{\norm{g_i(t)}^2}{\Filt_t}
 \le{} & 4(d\!+\!1)L_1^2\!\left(N\!\norm{\xbar(t)\!-\!x^\star}^2 \!+\! \norm{\bx_\perp}^2\right) \\
& +8N(d\!+\!1)L_1\mathcal H_f
 \!+\! 2N(d\!+\!2)L_1^2\beta_t^2.
\end{aligned}
\end{equation}
We also need to bound $\E\Econd{\norm{\gbar(t)}^2}{\Filt_t}$. By conditional independence of each $g_i(t)$ given $\Filt_t$ and $\E\Econd{g_i(t)}{\Filt_t}=\nabla f^{\beta_t}_{i}(x_i(t))$, where $f^{\beta_t}_{i}$ denotes the smoothed version of $f_i$ with smoothing radius $\beta_t$, we get
\begin{equation}
\label{eq:second_moment_gbar_sc_temp1}
\begin{aligned}
\E\Econd*{\norm{\gbar(t)}^2}{\Filt_t}
\!={} & \!\norm{\bar{m}(t)}^2 +
\mfrac{1}{N^2}\sum\nolimits_i \E\Econd*{\norm{\xi_i(t)}^2}{\Filt_t},
\end{aligned}
\end{equation}
where we denote $\bar{m}(t) = \E\Econd*{\gbar(t)}{\Filt_t}= \frac{1}{N}\sum_i \nabla f^{\beta_t}_{i}(x_i(t))$ and $\xi_i(t) = g_i(t)-\nabla f^{\beta_t}_{i}(x_i(t))$. We then have
\[
\begin{aligned}
& \norm{\bar{m}(t)-\nabla f(\xbar(t))} \\
\le{} & \mfrac{1}{N}\!\sum\nolimits_i\!\norm{\nabla\! f_i^{\beta_t}\!(x_i(t))
\!-\!\nabla\! f_i(x_i(t)) \!+\! \nabla\! f_i(x_i(t)) \!-\! \nabla\! f_i(\xbar(t))} \\
\le{} & \beta_t L_1 + \mfrac{L_1}{N}\sum\nolimits_i\|x_i(t)-\xbar(t)\|,
\end{aligned}
\]
where we used Lemma~\ref{lem:ball}(iv) and $L_1$-smoothness of $f_i$. Thus
\begin{equation}
\label{eq:offline-avg_gradient-bias_sq}
\norm{\bar{m}(t)-\nabla f(\xbar(t))}^2\le 2\beta_t^2L_1^2 + 2L_1^2N^{-1}\norm{\bx_\perp}^2,
\end{equation}
and
\[
\begin{aligned}
\norm{\bar{m}(t)}^2\le{} & 2\norm{\nabla f(\xbar(t))}^2 + 4\beta_t^2L_1^2 + 4L_1^2N^{-1}\norm{\bx_\perp}^2 \\
\le{} & 2L_1^2 \norm{\xbar(t)-x^\star}^2+ 4L_1^2\big(\beta_t^2 + N^{-1}\norm{\bx_\perp}^2\big).
\end{aligned}
\]
On the other hand, by the definition of $\xi_i(t)$, we have $\E\Econd{\norm{\xi_i(t)}^2}{\Filt_t}\le \E\Econd{\norm{g_i(t)}^2}{\Filt_t}$. Plugging the previous bounds into~\eqref{eq:second_moment_gbar_sc_temp1} and using~\eqref{eq:offline-moment} gives
\begin{equation}
\label{eq:gbar_second_moment_sc}
\begin{aligned}
\E\Econd*{\norm{\gbar(t)}^2}{\mkern-2mu \Filt_t}
\! \le{} & \!\Big(2\!+\!\mfrac{4(d\!+\!1)}{N}\Big) L_1^2\norm{\xbar(t) \!-\! x^\star}^2
\!+\! \mfrac{8(d\!+\!1)L_1\mathcal{H}_f}{N} \\
& \mkern-18mu
+4L_1^2\Big(1\!+\!\mfrac{d\!+\!1}{N}\Big)\mfrac{\norm{\bx_\perp}^2}{N}
\!+\!\beta_t^2 L_1^2\Big(4\!+\!\mfrac{2(d\!+\!2)}{N}\Big).
\end{aligned}
\end{equation}

\smallskip\noindent
\textbf{Step II: Recursions.} Combining~\eqref{eq:offline-avg_gradient-bias_sq} with 
$\inner{a}{b}\le(\mu/4)\norm{a}^2+\norm{b}^2/\mu$ gives
\[
\begin{aligned}
& \inner{\bar{m}(t)-\nabla f(\xbar(t))}{\! \xbar(t) - x^\star} \\
\ge{} & -\mfrac{\mu}{4}\!\norm{\xbar(t) \!-\! x^\star}^2
-\mfrac{2\beta_t^2L_1^2}{\mu}
-\mfrac{2L_1^2}{\mu N}\!\norm{\bx_\perp(t)}^2.
\end{aligned}
\]
Because $f$ is $\mu$-strongly convex and $\grad f(x^\star)=0$, we have
$\inner{\grad f(\xbar(t))}{\xbar(t)-x^\star}\ge\mu\norm{\xbar(t)-x^\star}^2$. Summing the previous two inner products leads to
\begin{align}
& \inner{\E\Econd*{\gbar(t)}{\Filt_t}}{\xbar(t) \!-\! x^\star}
\notag\\
\ge{} & \mfrac{3\mu}{4}\norm{\xbar(t)-x^\star}^2
-\mfrac{2\beta_t^2L_1^2}{\mu}
-\mfrac{2L_1^2}{\mu N}\norm{\bx_\perp(t)}^2.
\label{eq:offline-descent-direction}
\end{align}
Expanding $\norm{\xbar(t+1)-x^\star}^2$ along the average dynamics
$\xbar(t+1)=\xbar(t)-\eta_t\gbar(t)$ and substituting
\eqref{eq:gbar_second_moment_sc} and \eqref{eq:offline-descent-direction} yields
\begin{align}
 &\E\Econd*{\norm{\xbar(t+1)-x^\star}^2}{\Filt_t}\notag\\
={} & \! \norm{\xbar(t)\!-\!x^\star}^2 \!+\! \eta_t^2\E\Econd*{\norm{\gbar(t)}^2\!}{\!\Filt_t} \!-\! 2\eta_t\inner{\E\Econd*{\gbar(t)\!}{\!\Filt_t}}{\!\xbar(t)\!-\!x^\star} \notag \\
 \le{} & 
 \Big[1-\mfrac{3\mu\eta_t}{2}+ 2\eta_t^2L_1^2\!+\!\eta_t^2L_1^2\mfrac{4(d\!+\!1)}{N}\Big]
 \norm{\xbar(t)-x^\star}^2\notag\\
 &+
 4\eta_t L_1^2\left[\mfrac{1}{\mu}+\eta_t\Big(1\!+\!\mfrac{d\!+\!1}{N}\Big)\right]
 \mfrac{\norm{\xperp(t)}^2}{N}\notag\\
 &+
 \eta_t^2\mfrac{8(d \!+\! 1)L_1\mathcal H_f
\!+\! 2(d\!+\!2)\beta_t^2 L_1^2}{N}
 \!+\! \eta_t \beta_t^2 L_1^2 \left(
\mfrac{4}{\mu}\!+\!4\eta_t
\right).
\label{eq:offline-average-recursion}
\end{align}
For the consensus error,
$\xperp(t+1)=((P-J)\otimes I_d)\xperp(t)-\eta_t(M\otimes I_d)\bg(t)$.
Applying $\norm{p+q}^2\le(1+\alpha)\norm p^2+(1+\alpha^{-1})\norm q^2$ with
$\alpha=(1-\sigma^2)/(2\sigma^2)$ (the case $\sigma=0$ is immediate),
$\norm{M}_2\le1$, and \eqref{eq:offline-moment} gives
\begin{align}
 & \mfrac1N\E\Econd*{\norm{\xperp(t\!+\!1)}^2}{\Filt_t}
\notag\\
 \le{} &\mfrac{1\!+\!\sigma^2}{2}\!\mfrac{\norm{\xperp(t)}^2}{N}
\!+\! 4\eta_t^2\mfrac{1\!+\!\sigma^2}{1\!-\!\sigma^2}(d\!+\!1)L_1^2\Big[\!\norm{\xbar(t)\!-\!x^\star}^2\!+\!\mfrac{\norm{\bx_\perp}^2}{N}\Big]
\notag\\
 &+\mfrac{1\!+\!\sigma^2}{1\!-\!\sigma^2}\eta_t^2
 \left[8(d\!+\!1)L_1\mathcal H_f
+2(d\!+\!2)\beta_t^2L_1^2\right].
\label{eq:offline-consensus-recursion}
\end{align}

\smallskip\noindent
\textbf{Step III: Contraction and induction.} Temporarily denote $c=24L_1^2/(\mu(1-\sigma^2))$. Note that the condition on $\eta_t$ implies
\begin{equation}
\label{eq:step-size_bounds_sc}
\eta_t\le\min\!\left\{
\mfrac{\mu}{4L_1^2},\mfrac{\mu N}{16(d\!+\!1)L_1^2},\Big[\mfrac{\mu(1\!-\!\sigma^2)}{16c(1\!+\!\sigma^2)(d\!+\!1)L_1^2}\Big]^{\!\frac{1}{2}},\mfrac{1\!-\!\sigma^2}{4\mu}
\right\}.
\end{equation}
We now add \eqref{eq:offline-average-recursion} and
\eqref{eq:offline-consensus-recursion} multiplied by $c\eta_{t+1}$ to bound $\E\Econd{V_{t+1}}{\Filt_t}$. The coefficient of
$\norm{\xbar(t)-x^\star}^2$ in the sum is
\[
1-\mfrac{3\mu\eta_t}{2}+ 2\eta_t^2L_1^2\!+\!\mfrac{4(d\!+\!1)}{N}\eta_t^2L_1^2 \!+\! 4c\mfrac{1\!+\!\sigma^2}{1\!-\!\sigma^2}(d\!+\!1)L_1^2\eta_t^2\eta_{t+1}.
\]
By~\eqref{eq:step-size_bounds_sc} and $\eta_{t+1}\leq\eta_t$, we have $2\eta_t^2L_1^2\le \mu\eta_t/2$, $\frac{4(d+1)}{N}\eta_t^2L_1^2\le \mu \eta_t/4$, and $4c\mfrac{1\!+\!\sigma^2}{1\!-\!\sigma^2}(d\!+\!1)L_1^2\eta_t^2\eta_{t+1}\le \mu\eta_t/4$. Thus the coefficient of
$\norm{\xbar(t)-x^\star}^2$ can be bounded by $1-\mu\eta_t/2$, which is nonnegative since $\mu\le L_1$ and~\eqref{eq:step-size_bounds_sc} implies $\eta_t\le 1/(4\mu)$. The coefficient of $\norm{\xperp(t)}^2/N$ in the sum is
\[
\begin{aligned}
& c\eta_{t+1}\mfrac{1\!+\!\sigma^2}{2}
+ \mfrac{4\eta_t L_1^2}{\mu}+4\eta_t^2L_1^2+4\eta_t^2L_1^2\mfrac{d\!+\!1}{N} \\
& \qquad +4c\eta_t^2\eta_{t+1}\mfrac{1\!+\!\sigma^2}{1\!-\!\sigma^2}(d\!+\!1)L_1^2.
\end{aligned}
\]
Using the value of $c$ and $\eta_{t+1}\le \eta_t$, we get $c\eta_{t+1}(1+\sigma^2)/2\le c\eta_t - 12\eta_t L_1^2/\mu$. By~\eqref{eq:step-size_bounds_sc}, we have
\[
4\eta_t^2 L_1^2+4\eta_t^2L_1^2\mfrac{d\!+\!1}{N} + 4c\eta_t^2\eta_{t+1}\mfrac{1\!+\!\sigma^2}{1\!-\!\sigma^2}(d+1)L_1^2\le \mfrac{3}{2}\mu\eta_t,
\]
which is further bounded by $3\eta_t L_1^2/(2\mu)$.
Thus the coefficient of $\norm{\xperp(t)}^2/N$ can be bounded by
$c\eta_t - 13\eta_t L_1^2/(2\mu)\le c\eta_t(1-\mu\eta_t/2)$, where the last inequality can be justified by $\eta_t\le(1-\sigma^2)/(4\mu)$ from~\eqref{eq:step-size_bounds_sc}.
Eq.~\eqref{eq:step-size_bounds_sc} also implies $4/\mu+4\eta_t\le 5/\mu$. Hence
\begin{align*}
 \E\Econd*{V_{t+1}}{\Filt_t}
 \le{}&
 \Big(1-\mfrac{\mu}{2}\eta_t\Big)V_t
 +5\eta_t\beta_t^2 L_1^2/\mu
\notag\\
 & \mkern-54mu
+\eta_t^2[8(d\!+\!1)L_1\mathcal{H}_f \!+\! 2(d\!+\!2)L_1^2\beta_t^2]
\Big[
\mfrac{1}{N} \!+\! c\eta_t\mfrac{1\!+\!\sigma^2}{1\!-\!\sigma^2}
\Big].
\end{align*}
Plugging in $\eta_t=8/(\mu(t+b))$, $\beta_t\le \beta_{\mathrm{off}}/(t+b)$, taking total expectation, and using $t+b\ge b$, we can get from the above inequality that
\[
\E[V_{t+1}]
\le \Big(1-\mfrac{4}{t\!+\!b}\Big)\E[V_{t}]
+\mfrac{3C_1}{(t\!+\!b)^2} + \mfrac{2C_2}{(t\!+\!b)^3},
\]
where $C_1=512(d+1)L_1\mathcal{H}_f/(3\mu^2 N)$ and
\[
\begin{aligned}
C_2 ={} & \mfrac{20L_1^2\beta_{\mathrm{off}}^2}{\mu^2}
+(d\!+\!1)c\mfrac{1\!+\!\sigma^2}{1\!-\!\sigma^2}\mfrac{2048L_1\mathcal{H}_f}{\mu^3} + \mfrac{64(d\!+\!2)L_1^2\beta_{\mathrm{off}}^2}{N\mu^2 b} \\
& + \beta_{\mathrm{off}}^2 (d\!+\!2)c\mfrac{1\!+\!\sigma^2}{1\!-\!\sigma^2}\mfrac{512 L_1^2}{\mu^3 b^2}.
\end{aligned}
\]
By mathematical induction using $(t+1+b)^{-1}\ge (1-4/(t+b))(t+b)^{-1}+3(t+b)^{-2}$, $(t+1+b)^{-2}\ge (1-4/(t+b))(t+b)^{-2}+2(t+b)^{-3}$ and $(\frac{t+b}{t+1+b})^4\ge 1-4/(t+b)$, we obtain the bound~\eqref{eq:offline-V-rate}. Plugging the lower bound of $b$ into the definition of $C_2$ then leads to $C_2=O\big(
\frac{dL_1^3\mathcal{H}_f}{\mu^4(1-\sigma^2)^2}
+\beta_{\mathrm{off}}^2L_1^2/\mu^2\big)$.
The rate on $\E[\norm{\xbar(t)\!-\!x^\star}^2]$ is evident by $\norm{\xbar(t)\!-\!x^\star}^2\le V_t$.

To bound the consensus error, denote $e(t)\coloneq\frac1N\E\norm{\xperp(t)}^2$. Taking total expectations in
\eqref{eq:offline-consensus-recursion}, using~\eqref{eq:offline-V-rate} and $t+b\ge b$, and plugging in the choices of $\eta_t$ and $\beta_t$ gives
\[
e(t+1)\le \Big(\mfrac{1\!+\!\sigma^2}{2} + 4\eta_t^2\mfrac{1\!+\!\sigma^2}{1\!-\!\sigma^2}(d\!+\!1)L_1^2\Big)e(t)+ \mfrac{B}{(t+b)^2},
\]
where $B=O\big(\frac{dL_1(\mathcal{H}_f+L_1V_0)}{\mu^2(1-\sigma^2)}+(1-\sigma^2)\beta_{\mathrm{off}}^2\big)$ is a constant. Note that the choice of $\eta_t$ implies
\[
4\eta_t^2\mfrac{1\!+\!\sigma^2}{1\!-\!\sigma^2}(d\!+\!1)L_1^2
\le 
\frac{256}{\mu^2 b^2}\mfrac{1\!+\!\sigma^2}{1\!-\!\sigma^2}(d\!+\!1)L_1^2
\le \frac{1-\sigma^2}{4},
\]
and thus $e(t+1)\le \frac{3+\sigma^2}{4}e(t)+ \frac{B}{(t+b)^2}$. Furthermore, we also have $(t+1+b)^2/(t+b)^2\le (b+1)^2/b^2\le 8/(7+\sigma^2)$ from our choice of $b$. Based on these results, we can then show by induction that
\[
e(t)\le e(0)\Big(\mfrac{3+\sigma^2}{4}\Big)^{t} + \mfrac{8 B}{(1-\sigma^2)(t+b)^2}.
\]
Introducing $C_{\mathrm{con}}\coloneq 8B/(1-\sigma^2)$ completes the proof.

\section{Proofs for Results in Section~\ref{sec:privacy}}
\label{app:privacy}
 
\subsection*{Proof of Lemma~\ref{lem:density}}
Denote $G=(u^{\!\top}g)u$. Write $u=\varrho v$ with $\varrho =\norm{u}\in(0,1]$ and $v\in\Sphere^{d-1}$; under
$\Unif(\Ball^d)$, $\varrho$ and $v$ are independent, $v\sim\Unif(\Sphere^{d-1})$, and for any bounded Borel test function $h$,
\[
  \E[h(u)]
  =\frac{1}{\vol(\Ball^d)}
   \int_{\Sphere^{d-1}}\!\int_0^1 h(r\upsilon)\,r^{d-1}\,\ud r\,\ud\Omega(\upsilon),
\]
with $\Omega$ the surface measure. Since
$G=\varrho^2(v^{\!\top}g)\,v$ takes the same value at $(\varrho,v)$ and $(\varrho,-v)$, and
$\{v^{\!\top}g=0\}$ is $\Omega$-null, we may restrict to the hemisphere
$\{w:\ w^{\!\top}g>0\}$ and double:
\[
\begin{aligned}
 \E[h(G)]&=\frac2{\vol(\Ball^d)}\int_{\{w^{\!\top}g>0\}}\!\int_0^1
 h(r^2\alpha w)r^{d-1}\,\ud r\,\ud\Omega(w),
\end{aligned}
\]
where we denote $\alpha = w^{\!\top}g>0$. For fixed $w$, substitute $s=r^2\alpha\in(0,\alpha]$, so
$r=(s/\alpha)^{1/2}$ and $\ud r=\ud s/(2\sqrt{s\alpha})$:
\[
\begin{aligned}
 & \int_0^1h(r^2\alpha w)r^{d-1}\,\ud r
 =\int_0^\alpha h(sw)\Bigl(\frac s\alpha\Bigr)^{\frac{d-1}{2}}
 \frac{\ud s}{2\sqrt{s\alpha}} \\
 ={} & \frac12\int_0^\alpha h(sw)s^{\frac{d-2}{2}}\alpha^{-\frac d2}\,\ud s
 =\frac12\int_0^\alpha h(sw)(s\alpha)^{-\frac d2}s^{d-1}\,\ud s.
\end{aligned}
\]
Substituting this into the hemisphere integral, we 
recognize $s^{d-1}\,\ud s\,\ud\Omega(w)$ as the polar decomposition of
Lebesgue measure with respect to the variable $y=sw$. Then, using $s\alpha=g^{\!\top}y$ and observing that
$0<s\le\alpha$ is equivalent to $g^{\!\top}y\ge\norm y^2$, we obtain, for every
bounded Borel test function $h$,
\[
 \E[h(G)]=\int_{\R^d}h(y)
 \frac{\mathbf 1_{S_g\setminus\{0\}}(y)}{\vol(\Ball^d)}
 (g^{\!\top}y)^{-d/2}\,\ud y.
\]
This proves \eqref{eq:density}. The assigned value $\phi_g(0)=0$ is immaterial since a Lebesgue density is specified only almost everywhere.
 
\subsection*{Proof of Lemma~\ref{lem:cauchy}}

Denote $m=(g+g')/2$, $\Delta=g'-g$, and $\chi=\Dt(g,g')=\norm{\Delta}/\norm{m}$.

\smallskip
\noindent{}%
\emph{(i)} Since $g\ne g'$ we have $\Delta\ne0$. Choose orthonormal
vectors $e_1,e_2$ with $m=\norm{m}\,e_1$ and
$\Delta=\norm{\Delta}(\cos\theta\,e_1+\sin\theta\,e_2)$, where
$\cos\theta=\inner{\Delta}{m}/(\norm{\Delta}\norm{m})$ and,
by the choice of $e_2$, $\sin\theta\ge0$. Then,
\[
\begin{aligned}
 W=\frac{\Delta^{\!\top}v}{m^{\!\top}v}
 =\frac{\norm\Delta}{\norm m}\Bigl(\cos\theta+\sin\theta\frac{e_2^\top v}{e_1^\top v}\Bigr)
 =\chi(\cos\theta \!+\! \sin\theta\,Z),
\end{aligned}
\]
where we let $Z=e_2^\top v/e_1^\top v$.
Representing $v=\xi/\norm{\xi}$ with $\xi\sim\mathcal N(0,I_d)$ gives
$Z=e_2^\top\xi/e_1^\top\xi$, a ratio of independent standard normals, which is
standard Cauchy. An affine image of a standard Cauchy is Cauchy, so
$W\sim\mathrm{Cauchy}(\chi\cos\theta,\chi\sin\theta)
 =\mathrm{Cauchy}(\mut,\gat)$ with $\mut^2+\gat^2=\chi^2$ (when
$\sin\theta=0$ the law degenerates to the point mass at $\mut$).
Independence of $\varrho$ and $W$ holds because $W$ is a function of $v$ alone while
$\varrho$ and $v$ are independent.
 
\medskip\noindent
\emph{(ii)} Since $m\ne0$ and $g\ne0$, the event
$\{m^{\!\top}v\ne0,\ g^{\!\top}v\ne0\}$ has probability one. On this
event, using $g=m-\Delta/2$ and $g'=m+\Delta/2$,
\[
  X=\frac{g'^{\top}v}{g^{\!\top}v}
   =\frac{m^{\!\top}v+\tfrac12\Delta^{\!\top}v}
         {m^{\!\top}v-\tfrac12\Delta^{\!\top}v}
   =\frac{1+W/2}{1-W/2}
   =\frac{2+W}{2-W}.
\]
On $\{W<2\}$ the map $W\mapsto\frac{2+W}{2-W}$ is increasing with range
$(-1,\infty)$, and on $\{W>2\}$ it takes values in $(-\infty,-1)$.
Hence, by $e^{2\epsilon/d}\ge1$, $X>e^{2\epsilon/d}$ forces $W<2$ and is then equivalent to $2+W>e^{2\epsilon/d}(2-W)$, i.e.,
$W>2(e^{2\epsilon/d}-1)/(e^{2\epsilon/d}+1)=2\tanh(\epsilon/d)$, giving
$\{X>e^{2\epsilon/d}\}=\{2\tanh(\epsilon/d)<W<2\}$. Similarly, $X\le0$ holds iff
$W\le-2$ or $W>2$, i.e.\ $\{|W|\ge2\}$ up to null sets; and $0<X<1$ iff
$-2<W<0$.

\subsection*{Proof of Theorem~\ref{thm:single}}
If $g=g'$, the two laws coincide and the result is immediate. Otherwise,
$\Dt(g,g')<2$ rules out either gradient being zero and also gives $g+g'\ne0$,
so Lemmas~\ref{lem:density} and~\ref{lem:cauchy} apply. For brevity, write
\[
 D_\epsilon=\{y\in\R^d:\phi_g(y)>e^\epsilon\phi_{g'}(y)\}.
\]
Outside $D_\epsilon$ the desired pointwise density inequality holds by
construction. To bound $\mathsf P_g(D_\epsilon)$, let
$Y=(u^{\!\top}g)u\sim\mathsf P_g$ and
write $u=\varrho v$, $X=g'^{\top}v/g^{\!\top}v$ as in Lemma~\ref{lem:cauchy}.
Then $Y\in S_g$ almost surely, and
\[
\begin{aligned}
 Y\in S_{g'}&\iff g'^{\top}Y\ge\norm{Y}^2\\
 &\iff \varrho^2(g^{\!\top}v)(g'^{\top}v)\ge \varrho^4(g^{\!\top}v)^2\\
 &\iff X\ge \varrho^2.
\end{aligned}
\]
On $\{Y\in S_g\cap S_{g'}\}$ the ratio \eqref{eq:llr} equals $X^{d/2}$,
which exceeds $e^{\epsilon}$ iff $X>e^{2\epsilon/d}$. Hence
\[
\begin{aligned}
 \{Y\in D_\epsilon\}
 \subseteq{} & \{X\le0\}\cup\{0<X<1,\, \varrho^2>X\}\\
 &\cup\{X>e^{2\epsilon/d}\},
\end{aligned}
\]
where we used $\varrho\le1$ to place $\{0<X<\varrho^2\}$ inside the middle event, and
noted that on $\{X>e^{2\epsilon/d}\}$ the support condition $X\ge \varrho^2$ is
automatic since $X>1\ge \varrho^2$. By Lemma~\ref{lem:cauchy}(ii) these events
are $\{|W|\ge2\}$, $\{-2<W<0,\varrho^2>X\}$, and
$\{\tau_\epsilon<W<2\}\subseteq\{W>\tau_\epsilon\}$. We bound the probability of each event in
turn; recall from Lemma~\ref{lem:cauchy}(i) that
$W\sim\mathrm{Cauchy}(\mut,\gat)$ with $\mut^2+\gat^2=\Dt(g,g')^2$, and that
$W$ is independent of $\varrho$.
The calculations below are written for $\gat>0$; the degenerate case
$\gat=0$ follows by a direct point-mass check, equivalently by continuity
as $\gat\downarrow0$. For notational simplicity, we denote $\tau_\epsilon = 2\tanh(\epsilon/d)$ and $\chi=\Dt(g,g')$.
 
\medskip\noindent
\textbf{The high-ratio event $\{W>2\tanh(\epsilon/d)\}$.}
Since $\mut\le \chi<\tau_\epsilon$ (by hypothesis), we have
$\tau_\epsilon-\mut>0$, and the
Cauchy right tail has the closed form
\[
  \Prob(W>\tau_\epsilon)
  =\frac12-\frac1\pi\arctan\frac{\tau_\epsilon-\mut}{\gat}
  =\frac1\pi\arctan\frac{\gat}{\,\tau_\epsilon-\mut\,}.
\]
Maximize
$\gat/(\tau_\epsilon-\mut)
 =\chi\sin\theta/(\tau_\epsilon-\chi\cos\theta)$ over
$\theta\in[0,\pi]$ on the circle $\mut^2+\gat^2=\chi^2$: The derivative of
the objective is proportional to $\tau_\epsilon\cos\theta-\chi$, so the
maximum is at $\cos\theta=\chi/\tau_\epsilon$, with value
$\chi/\sqrt{\tau_\epsilon^2-\chi^2}$. Using
$\arctan\bigl(\chi/\sqrt{\tau_\epsilon^2-\chi^2}\bigr)
 =\arcsin(\chi/\tau_\epsilon)$, we get
\[
  \Prob(W>\tau_\epsilon)
  \le\frac1\pi\arcsin\frac{\chi}{\tau_\epsilon}.
\]
 
\medskip\noindent
\textbf{The directional-escape event $\{|W|\ge2\}$.}
Since $|\mut|\le\chi<2$, both $2\pm\mut>0$, and by the Cauchy tail
bound $\arctan z\le z$,
\[
\begin{aligned}
\Prob(W\le-2)+\Prob(W\ge2)
 &\le\frac1\pi\Bigl(\frac\gat{2+\mut}+\frac\gat{2-\mut}\Bigr)\\
 &=\frac\gat\pi\frac4{4-\mut^2}\le\frac{4\chi/\pi}{4-\chi^2},
\end{aligned}
\]
using $\gat\le\chi$ and $\mut^2\le\chi^2$.
 
\medskip\noindent
\textbf{The radial-escape event $\{-2<W<0,\ \varrho^2>X\}$.}
Since $\varrho$ is independent of $W$
(Lemma~\ref{lem:cauchy}(i)) with $\Prob(\varrho\le s)=s^d$ for $s\in[0,1]$, on the event
$\{W\in(-2,0)\}=\{X\in(0,1)\}$ we have
$\Prob\Pcond{\varrho^2>X}{W}=1-X^{d/2}$. Then, by
$1-X^{d/2}\le\tfrac d2(1-X)$ and
$1-X=\tfrac{-2W}{2-W}\le-W$ on $\{W\in(-2,0)\}$, we get
\[
  \Prob\bigl(-2<W<0,\,\varrho^2>X\bigr)
  \le\frac d2\,\E\bigl[(-W)\,\mathsf{1}_{\{-2<W<0\}}\bigr],
\]
where $\mathsf{1}_{\{-2<W<0\}}$ denotes the corresponding indicator function. Next, set $W'=-W\sim\mathrm{Cauchy}(-\mut,\gat)$. Direct integration of the
Cauchy density gives
\[
\begin{aligned}
 \E\!\left[W'\mathsf{1}_{\{0<W'<2\}}\right]
 ={} & -\mut\,\Prob(0<W'<2)\\
 &+\frac\gat{2\pi}\ln\frac{(2+\mut)^2+\gat^2}{\mut^2+\gat^2}.
\end{aligned}
\]
Using $\mut^2+\gat^2=\chi^2$,
$(2+\mut)^2+\gat^2=4+4\mut+\chi^2\le(2+\chi)^2$, together with
$-\mut\,\Prob(\cdot)\le|\mut|\le\chi$ and $\gat\le\chi$, we obtain
\[
  \E\bigl[(-W)\,\mathsf{1}_{\{-2<W<0\}}\bigr]
  \le\chi+\frac{\chi}{\pi}\ln\frac{2+\chi}{\chi}.
\]
Multiplying by $d/2$ gives the claim.

Summing the three bounds gives
$\mathsf P_g(D_\epsilon)\le\delta_{\mathrm{step}}(\epsilon;\chi)$. Therefore, for
every measurable $A$,
\[
\begin{aligned}
 \mathsf P_g(A)
 &=\int_{A\setminus D_\epsilon}\phi_g(y)\,dy
   +\mathsf P_g(A\cap D_\epsilon)\\
 &\le e^\epsilon\int_{A\setminus D_\epsilon}\phi_{g'}(y)\,dy
   +\delta_{\mathrm{step}}(\epsilon;\chi)\\
 &\le e^\epsilon\mathsf P_{g'}(A)
   +\delta_{\mathrm{step}}(\epsilon;\chi).
\end{aligned}
\]
Exchanging $g$ and $g'$ proves the second inequality in
\eqref{eq:single-round-dp}.
Each term in \eqref{eq:deltastep} is increasing in $\chi$ on $[0,2\tanh(\epsilon/d))$, and $\chi\log((2+\chi)/\chi)\to0$ as $\chi\downarrow0$, proving the
remaining claims.

\bibliographystyle{IEEEtran}
\bibliography{references}
\end{document}